\documentclass[12pt,psamsfonts,leqno,oneside,letterpaper]{amsart}
\usepackage[dvips,text={6.5truein,9truein},left=1truein,top=1truein]{geometry}
\usepackage{amssymb,amsmath,amscd,enumerate}
\usepackage{graphicx}
\usepackage{url}

\usepackage[colorlinks,linkcolor=blue,citecolor=blue,pdfstartview=FitH]{hyperref}
\input xy
\xyoption{all}
\SelectTips{cm}{12}

\usepackage{color}

\theoremstyle{definition}
\newtheorem{para}{}[section]

\newtheorem{subclaim}[equation]{}
\newtheorem{remark}[para]{Remark}

\newtheorem{definition}[para]{Definition}

\newtheorem{claim}[equation]{Claim}

\theoremstyle{plain}

\newtheorem{lemma}[para]{Lemma}
\newtheorem{proposition}[para]{Proposition}
\newtheorem{corollary}[para]{Corollary}

\numberwithin{equation}{para}
\numberwithin{figure}{section}

\newcommand\calAC{\mathcal{AC}}
\newcommand\calBC{\mathcal{BC}}
\newcommand\calHC{\mathcal{HC}}
\newcommand\calHI{\mathcal{HI}}
\newcommand\calc{{\mathcal C}}
\newcommand\calE{\mathcal{E}}

\newcommand\cali{\mathcal{I}}

\newcommand\calh{{\mathcal H}}

\newcommand\co{\colon\thinspace}
\newcommand\bd{\mathbf{d}}

\newcommand\dist{\mathit{dist}}

\begin{document}

\title{The geometry of cyclic hyperbolic polygons}

\author{Jason DeBlois}
\address{Department of Mathematics\\University of Pittsburgh}
\email{jdeblois@pitt.edu}
\thanks{This paper was mostly written during a period of partial NSF support}

\begin{abstract}  A hyperbolic polygon is defined to be \textit{cyclic}, \textit{horocyclic}, or \textit{equidistant} if its vertices lie on a metric circle, horocycle, or a component of the equidistant locus to a hyperbolic geodesic, respectively.  Convex such $n$-gons are parametrized by the subspaces of $(\mathbb{R}^+)^n$ that contain their side length collections, and area and circumcircle or ``collar'' radius determine symmetric, smooth functions on these spaces.  We give formulas for and bounds on the derivatives of these functions, and make some observations on their behavior.  Notably, the monotonicity properties of area and circumcircle radius exhibit qualitative differences on the collection of centered vs non-centered cyclic polygons, where a cyclic polygon is \textit{centered} if it contains the center of its circumcircle in its interior.\end{abstract}

\maketitle


An $n$-tuple $(d_0,\hdots,d_{n-1})$ of positive real numbers is the side length collection of a compact convex $n$-gon $P$ in the hyperbolic plane $\mathbb{H}^2$ that is cyclic or horocyclic, or has all vertices equidistant from a fixed geodesic, as long as $d_i \leq \sum_{j\neq i} d_j$ for each $i$. Moreover, such a polygon $P$ is uniquely prescribed up to orientation-preserving isometry of $\mathbb{H}^2$ by this $n$-tuple up to cyclic permutation.  It is cyclic if and only if for each $i$, 
$$\sinh(d_i/2) < \sum_{j\neq i} \sinh(d_j/2);$$ 
and $P$ is horocyclic if and only if for some $i$ equality holds above.  

This fact seems to have been independently rediscovered several times.  It was recorded by the late W.W.~Stothers \cite{Stothers}, then by J.M.~Schlenker (see \cite[p.~2175]{Schlenker}) as part of a larger project.  Schlenker's results on cyclic polygons were re-proved by R.~Walter \cite{Walter_poly1}, \cite{Walter_poly2}.  Propositions \ref{existence uniqueness}, \ref{horocyclic parameters} and \ref{equidistant parameters} here combine for a self-contained case-by-case proof.  Our approach is standard; compare eg.~Robbins \cite{Robbins} (in the Euclidean setting).  All aspects are elementary, but the cyclic case is somewhat subtle; see the discussion at the beginning of Section \ref{existence}.

Our main goal here is to make some observations on the qualitative behavior of area and circumcircle radius, thought of as functions on a subset $\calAC_n$ of $(\mathbb{R}^+)^n$ that parametrizes cyclic $n$-gons by side length, for each $n\geq 3$.  We also record relatively simple formulas for their derivatives, and give useful bounds in some cases.  In view of the facts above we define:\begin{align*}
   \mathcal{AC}_n & = \left\{(d_0,\hdots,d_{n-1}) \in (\mathbb{R}^+)^n\,|\, \sinh (d_i/2) < \sum_{j\neq i} \sinh(d_j/2) \ \mbox{for each}\ i \in \{0,\hdots,n-1\}\right\} \end{align*}
One of our key themes is that the functions on $\calAC_n$ determined by area and circumcircle radius  exhibit qualitatively different behaviors on \textit{centered} versus non-centered cyclic $n$-gons.  We say a cyclic polygon is centered if it contains its circumcircle center in its interior.  The collection $\calc_n\subset\calAC_n$ defined below parametrizes centered $n$-gons.  Here for $J\geq d/2$ let $\theta(d,J) = 2\sin^{-1}\left(\sinh(d_i/2)/\sinh J\right)$.\begin{align*}
  \mathcal{C}_n & = \left\{(d_0,\hdots,d_{n-1}) \in (\mathbb{R}^+)^n\,|\,\sum_{i=0}^{n-1}\theta(d,D/2) > 2\pi,\ \mbox{where}\ D=\max\{d_i\}_{i=0}^{n-1} \right\} \end{align*}
The qualitative differences we mentioned above are visible in the bounds below:

\newcommand\SmoothJ{For $n\geq3$, the function $J\co \mathcal{AC}_n\to\mathbb{R}^+$ that records circumcircle radius is smooth and symmetric.  For $\bd=(d_0,\hdots,d_{n-1})\in\calAC_n$:
$$ \left\{\begin{array}{ll}  0 < \frac{\partial J}{\partial d_i}(\bd) < 1/2 & \mbox{if}\ \bd\in\calc_n,\ \mbox{for any}\ i \\
  \frac{\partial J}{\partial d_{i_0}}(\bd) > 1/2 & \mbox{if}\ \bd\in\calAC_n - (\calc_n\cup\calBC_n)\ \mbox{and}\ d_{i_0} = \max\{d_i\}_{i=0}^{n-1}\\
  \frac{\partial J}{\partial d_j}(\bd) < 0 & \mbox{if}\ \bd\in\calAC_n - (\calc_n\cup\calBC_n)\ \mbox{and}\ d_j\neq\max\{d_{i}\}_{i=0}^{n-1}
  \end{array}\right.$$
Furthermore, if $d_i > d_j$ then $\left|\frac{\partial J}{\partial d_i}(\bd)\right| > \left|\frac{\partial J}{\partial d_j}(\bd)\right|$.}
\newtheorem*{SmoothJprop}{Proposition \ref{smooth J}}
\begin{SmoothJprop}\SmoothJ\end{SmoothJprop}

The subspace $\calBC_n$ of $\calAC_n$ referenced above is the frontier of $\calc_n$ in $\calAC_n$, a codimension-one submanfold $\calBC_n$ that parametrizes \textit{semicyclic} $n$-gons: those whose longest edge is a diameter of their circumcircle.  (See Propositions \ref{Cn and BCn} and \ref{submanifold}; the term ``semicyclic'' is due to Maley--Robbins--Roskies in the Euclidean setting \cite{MRR}.)  Values of the $\partial J/\partial d_i$ on $\calBC_n$ are thus determined by continuity and the formulas above.

We also give a Schl\"afli-type formula on areas of cyclic $n$-gons.  Whereas the classical Schl\"afli formula expresses the change in area of a first-order deformation of polygons in terms of angle variations, ours is in terms of side length.  And since cyclic polygons are parametrized by their side length collections, we can simply record partial derivatives.

\newcommand\SmoothD{For $n\geq3$, the function $D_0\co\calAC_n\to\mathbb{R}^+$ that records hyperbolic area is smooth and symmetric.  For $\bd= (d_0,\hdots,d_{n-1})\in\calAC_n$:
$$ \frac{\partial D_0}{\partial d_i}(\bd) = \left\{\begin{array}{ll}
  -\sqrt{\frac{1}{\cosh^2(d_i/2)} - \frac{1}{\cosh^2 J(\bd)}} & 
  \mbox{if}\ \bd\in\calAC_n-\calc_n\ \mbox{and}\ d_i=\max\{d_j\}_{j=0}^{n-1}\\ 
  \ \ \sqrt{\frac{1}{\cosh^2(d_i/2)} - \frac{1}{\cosh^2 J(\bd)}} & 
  \mbox{otherwise}\end{array}\right.$$}
\newtheorem*{SmoothDprop}{Proposition \ref{smooth D}}
\begin{SmoothDprop}\SmoothD\end{SmoothDprop}

Note that the circumcircle radius of a semicyclic $n$-gon is half its longest side length, so the two  formulas above are compatible with continuous partial derivatives of $D_0$; in particular, the first is $0$ everywhere on $\calBC_n$.  One consequence of Proposition \ref{smooth D} is that area of centered and semicyclic $n$-gons is monotonic in their side lengths:

\newcommand\Monotonicity{For $n \geq 3$ and $\bd = (d_0,\hdots,d_{n-1})$ and $\bd'=(d_0',\hdots,d_{n-1}')$ in $\calc_n\cup\calBC_n$, if after a permutation $d_{i} \leq d_i'$ for all $i$, and $d_i<d_i'$ for some $i$, then $D_0(\bd) < D_0(\bd')$.}
\newtheorem*{monotonicitycor}{Corollary \ref{monotonicity}}
\begin{monotonicitycor}\Monotonicity\end{monotonicitycor}

Section \ref{to infinity} extends the considerations above to horocyclic polygons in $\mathbb{H}^2$.  We parametrize horocyclic $n$-gons by a codimension-one submanifold $\calHC_n$ of $(0,\infty)^n$, which is the frontier of $\calAC_n$ there.  The function $J\co\calAC_n\to\mathbb{R}^+$ blows up approaching $\calHC_n$ (see Proposition \ref{radius up_n_down}), suggesting the geometric interpretation that horocyclic $n$-gons are limits of sequences of cylic $n$-gons whose circumcircle radii go to infinity.

Section \ref{beyond} considers equidistant polygons, parametrized by the set $\calE_n$ below which, like $\calAC_n$, has nonempty interior in $(0,\infty)^n$.  In fact the closures there of $\calE_n$ and $\calAC_n$ intersect in $\calHC_n$. 
$$\calE_n = \left\{ (d_0,\hdots,d_n)\in(0,\infty)^n\,|\, \sinh(d_{i_0}/2) > \sum_{i\neq i_0} \sinh(d_i/2)\ \mbox{but}\ d_{i_0} \leq \sum_{i\neq i_0} d_i,\ \mbox{for some}\ i_0 \right\}$$
There is a ``collar radius'' function $J\co\calE_n\to[0,\infty)$ analogous to the circumcircle radius function on $\calAC_n$, with similar behavior, see Proposition \ref{equi smooth J}.  The area function $D_0$ extends continuously to $\calAC_n\cup\calHC_n\cup\calE_n$, by Propositions \ref{horocyclic defect} and \ref{equi defe deri}. Proposition \ref{equi defe deri} also gives a Schl\"afli-type formula similar to the one above for $D_0$ on $\calE_n$.  Combining Propositions \ref{smooth D} and \ref{equi defe deri} with Schlenker's results yields:

\newcommand\MaxArea{For fixed positive real numbers $d_1,\hdots,d_{n-1}$, among all hyperbolic $n$-gons with $n-1$ sides of these lengths, area is maximized by the semicyclic ones with final side length greater than $\max\{d_i\}$.}
\newtheorem*{maxareacor}{Corollary \ref{max area}}
\begin{maxareacor}\MaxArea\end{maxareacor}

The analogous result for Euclidean polygons is mentioned by Maley--Robbins--Roskies on p.~672 of \cite{MRR}, but we do not know a reference for a proof.

Section \ref{to infinity} also describes \textit{horocyclic ideal} polygons, which have all vertices on a horocycle except one at its ``ideal point'' (see Definition \ref{horocyclic order}).  These are natural limits for families of cyclic polygons with certain edge lengths going to infinity, and the area function $D_0$ extends continuously to their parameter space $\calHI_n$.  See Propositions \ref{horocyclic defect} and \ref{horocyclic decomp}.

Finally, in Section \ref{degenerations} we describe how cyclic $n$-gons degenerate to cyclic $m$-gons (for $m<n$) as some side lengths approach zero.  In terms of our parameter spaces, the frontier of $\calAC_n\cup\calHC_n$ in $[0,\infty)^n$ is a union of copies of $\calAC_m\cup\calHC_m$, for $3\leq m<n$, and a degenerate part (see Lemma \ref{shrinking limit}).  The main results are that the radius and area functions on $\calAC_n$ extend continuously to those on $\calAC_m$, see Lemma \ref{shrinking J} and Corollary \ref{shrinking D} respectively.

The results here are an important tool in \cite{DeB_Voronoi}, on Delaunay tessellations of hyperbolic surfaces.  The $2$-cells of such a decomposition are cyclic, horocyclic or equidistant polygons (see \cite{DeB_Delaunay}).

In the Euclidean setting, a body of work on giving explicit formulas for area and circumcircle radius reaches back to classical times.  Heron's formula, recorded circa A.D.~60, and Brahmagupta's formula from the seventh century give the areas of cyclic Euclidean triangles and quadrilaterals, respectively, as functions of their side lengths.  In 1828, A.F.~M\"obius gave formulas for circumcircle radius \cite{Mobius}.  More recently, D.~Robbins made a series of conjectures on a polynomial relation extending the Heron and Brahmagupta formulas to arbitrary $n$-gons \cite{Robbins}, which were subsequently proved in independent work of several authors (see eg.~\cite{Pak}).

A hyperbolic version of Heron's formula was only recorded in 1969 by S.~Bilinski \cite{Bilinski}, then independently rediscovered by Stothers.  Even more recently, A.~Mednykh found a hyperbolic Brahmagupta formula \cite{Mednykh}.  We know of no more literature on cyclic hyperbolic polygons beyond what has been cited above.  Conversely, we are not aware of recorded Euclidean analogs to our results (though these should not be difficult).

\subsection*{Acknowledgement} The author thanks the anonymous referee for his or her careful reading and suggestions that have improved the paper.

\section{Calculus (and some geometry)}\label{existence}

A straightforward ``method of continuity'' approach to proving existence of cyclic polygons with given side lengths is neatly sketched by D.~Robbins in the Euclidean setting \cite[p.~526]{Robbins}:

\quad\parbox[t]{445pt}{\it
Imagine a circle of variable radius and let us try to inscribe a polygon with sides of the given lengths in the circle by picking an arbitrary starting point and laying out the edges, one at a time, with the given lengths.  When the radius is too large, we will not reach the starting point when we have used up all the sides.  As we decrease the radius there will come a time when we return exactly to our starting point.}

\smallskip
This only works given two assumptions, though (both pointed out by Robbins): first, that all side lengths are nearly equal; and second, that all edges go in the same direction. That is, each edge joins its initial vertex $x_i$ to the possibility for $x_{i+1}$ that is closer on the circle to $x_i$ in the direction of some fixed orientation.

If one drops the first assumption but keeps the second then sometimes the circle radius can decrease as far as possible without the starting point being reached.  But on dropping the second assumption uniqueness can fail, since if the circle radius is $J$ and $d<2J$ then there are two possible edges of length $d$ with any given initial point.

However for \textit{convex} polygons it is not hard to see that all edges except possibly the longest must go in the same direction, and Robbins' strategy can be successfully modified as follows:

\quad\parbox[t]{445pt}{\it
Having laid out all edges in the same direction, if on decreasing the radius as far as we can (to half the length of the longest side) we have not yet returned to our starting point then re-inflate the circle, now forcing the longest edge to go in the opposite direction from the others.  As we increase the radius there will come a time when we return exactly to our starting point (assuming the edge lengths satisfy the proper inequality).}

\begin{figure}
\includegraphics{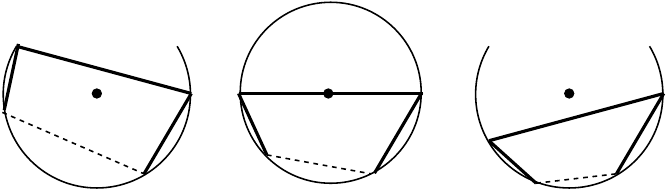}
\caption{From left to right, circle radius decreases then increases, but distance between endpoints of the union of edges continues to decrease.}
\label{n vs s}
\end{figure}

\smallskip
See Figure \ref{n vs s}.  Below we prove existence and uniqueness of convex cyclic hyperbolic $n$-gons by formalizing a slight variation on this approach: rather than laying out all $n$ edges we lay out $n-1$, and we vary the circle radius until the distance between start and end points equals the length of the omitted edge.  (Doing it this way helps with the calculus.) 

Before we begin we note that in \cite{Robbins}, Robbins considers Euclidean cyclic $n$-gons with no convexity requirement, whereupon uniqueness does indeed fail; see Diagram 1 of \cite{Robbins}.  In fact the extent of its failure determines the degree of the ``generalized Heron polynomials'' that are the main objects of study in \cite{Robbins}.

\subsection{Existence and uniqueness}\label{uniqueness}  Let us begin by specifying a model of the hyperbolic plane and recalling a few basic facts, which can be found in eg.~\cite{Ratcliffe} (or many others).

\begin{definition}\label{uhp} The \textit{upper half-plane model} of the hyperbolic plane is the open subset $\mathbb{H}^2\doteq\{z\in\mathbb{C}\,|\, \Im z>0\}$ of $\mathbb{C}$, inheriting the standard orientation and endowed with the Riemannian metric $\frac{dx^2+dy^2}{y^2}$ (identifying $\mathbb{C}$ with $\mathbb{R}^2$).  For $x$ and $x'$ in $\mathbb{H}^2$ let $\dist(x,x')$ refer to the distance from $x$ to $x'$ in the resulting path metric.\end{definition}

The orientation-preserving isometry group of $\mathbb{H}^2$ is $\mathrm{PSL}_2(\mathbb{R})$, acting by M\"obius transformations.  This action is transitive on $\mathbb{H}^2$, hence also on hyperbolic circles of fixed radius.  Each hyperbolic circle is a Euclidean circle contained in $\mathbb{H}^2$, though its hyperbolic and Euclidean center and radius do not coincide.

\begin{definition}\label{cyclic order}  For points $x$ and $y$ on a circle $C$, the \textit{counterclockwise} arc $[x,y]$ from $x$ to $y$ is the arc of $C$ bounded by $x$ and $y$ that points from $x$ to $y$ in the counterclockwise orientation on $C$.  The \textit{angle} from $x$ to $y$ is the angle subtended by $[x,y]$ at the center of $C$.  We will say a collection $\{x_0,\hdots,x_{n-1}\}\subset C$ is \textit{cyclically ordered} if for all $i>0$, $(x_{i-1},x_i]$ contains no $x_j$ for $j\neq i$.

A \textit{(convex) cyclic $n$-gon with side length collection} $(d_0,\hdots,d_{n-1})\in(0,\infty)^n$ is a cyclically ordered collection $\{x_0,\hdots,x_{n-1}\}$ on a fixed hyperbolic circle, called the \textit{circumcircle}, such that $\dist(x_{i-1},x_i) = d_i$ for each $i>0$ and $\dist(x_0,x_{n-1}) = d_0$.  An \textit{isometry} from a cyclic $n$-gon $\{x_0',\hdots,x_{n-1}'\}$ to $\{x_0,\hdots,x_{n-1}\}$ is an isometry of $\mathbb{H}^2$ that takes $x_i'$ to $x_{\sigma(i)}$ for each $i$, where $\sigma$ is a fixed cyclic permutation.\end{definition}

\begin{remark}\label{important}  If an isometry takes $x_i'$ to $x_{\sigma(i)}$ for each $i$, where $\sigma$ is a fixed cyclic permutation, then also $d_i' = d_{\sigma(i)}$ for each $i$.  In particular, \textit{cyclically relabeling} a cyclic $n$-gon $\{x_0,\hdots,x_{n-1}\}$ by a cyclic permutation $\sigma$ produces a cyclic $n$-gon $\{x_{\sigma(0)},\hdots,x_{\sigma(n-1)}\}$, isometric to the original via the identity map, with side length collection that of the original relabeled by $\sigma$.\end{remark}

Hyperbolic trigonometry relates the distance between points on a circle to the angle from one to another and the circle radius.  Here is the key function:

\begin{lemma}\label{angle function}  For $d>0$ and $J\geq d/2$, let $T$ be a hyperbolic triangle with two sides of length $J$ and a third of length $d$. Its interior angle at the vertex opposite the side of length $d$ is:
$$\theta(d,J) = 2\sin^{-1}(\sinh(d/2)/\sinh J)\in [0,\pi]$$
This is a continuous function on $\{(d,J)\,|\,0<d\leq 2J\}$, smooth on its interior.  For fixed $d>0$, $\theta(d,J)$ decreases in $J$ on $[d/2,\infty)$, with $\theta(d,d/2) = \pi$, $\lim_{J\to\infty} \theta(d,J) = 0$ and
$$ \frac{\partial}{\partial J} \theta(d,J) 
  = \frac{-2\sinh(d/2)\cosh J}{\sinh J\sqrt{\sinh^2 J - \sinh^2(d/2)}} 
  = -2\coth J \tan(\theta(d,J)/2) $$\end{lemma}

\begin{proof} Applying the hyperbolic law of sines (see \cite[Th.~3.5.2]{Ratcliffe}) to a triangle obtained from $T$ by dropping a perpendicular to the side with length $d$ from the vertex opposite it yields:\begin{align}\label{LOS leader}  \sinh(d/2) = \sin (\theta(d,J)/2)\sinh J \end{align}
The formula for $\theta(d,J)$ follows. Continuity and smoothness of $\theta$ follow from the fact that the inverse sine is continuous on $[-1,1]$ and smooth on $(-1,1)$.  By direct computation, $\theta(d,d/2) = 2\sin^{-1}(1) = \pi$ and for fixed $d>0$, $\lim_{J\to\infty} \theta(d,J) = 2\sin^{-1}(0) = 0$.  That $\theta(d,J)$ decreases in $J$ for fixed $d>0$ follows from the partial derivative computation, which is straightforward.\end{proof}

The result below describes the two ways to put a cyclic $n$-gon in a circle of fixed radius, specifying all but one of its side lengths.

\begin{proposition}\label{omnibus1}  For any $d_1,\hdots,d_{n-1}\in(0,\infty)$ and $J\geq\max\{d_i/2\}$, and any hyperbolic circle $C$ of radius $J$, there exists $\{x_0,\hdots,x_{n-1}\}\subset C$ such that $\mathit{dist}(x_{i-1},x_i) = d_i$ for each $i>0$.   One such collection has the angle $\theta_i$ from $x_{i-1}$ to $x_i$ equal to $\theta(d_i,J)$ for each $i>0$.  For any fixed $i_0$, there is another such that $\theta_{i_0} = 2\pi-\theta(d_{i_0},J)$, but $\theta_i = \theta(d_i,J)$ for each $i\neq i_0$.  Each  of these collections is determined up to hyperbolic isometry by the choice of $\theta_i$.  The former is cyclically ordered if and only if $\sum_{i=1}^{n-1} \theta(d_{i},J) <2\pi$, the latter if and only if $\theta(d_{i_0},J) > \sum_{i\neq i_0}\theta(d_i,J)$.

If $\{x_0,\hdots,x_{n-1}\}$ is a convex cyclic $n$-gon with side length collection $(d_0,\hdots,d_{n-1})$ on a hyperbolic circle of radius $J$, then $J\geq \max\{d_i/2\}$ and for each $i$, the angle $\theta_i$  from $x_{i-1}$ to $x_i$ is either $\theta(d_i,J)$ or $2\pi-\theta(d_i,J)$.  If $\theta_i = \theta(d_i,J)$ for all $i >0$ then:
$$ d_0 = \ell^n(J,d_1,\hdots,d_{n-1}) \doteq 2\sinh^{-1}\left[\sinh J\sin\left(\frac{1}{2}\sum_{i=1}^{n-1}\theta(d_i,J)\right)\right] $$
If  $\theta_i = 2\pi - \theta(d_{i_0},J)$ for some $i_0>0$, then $\sum_{i\neq i_0}\theta(d_i,J) <\theta(d_{i_0},J)$, so in particular $d_{i_0} > d_i$ for all $i\neq i_0$.  In this case:
$$ d_0 = \ell^s(J,d_1,\hdots,d_{n-1})\doteq 2\sinh^{-1}\left[\sinh J\sin\left(\frac{1}{2}\theta(d_{i_0},J) - \frac{1}{2}\sum_{i\neq i_0} \theta(d_i,J)\right)\right] $$
\end{proposition}

The letters ``$n$'' and ``$s$'' decorating $\ell$ above stand for ``non-separating'' and ``separating'', respectively, and correspond to the left and right sides of Figure \ref{n vs s}.  On the right, the longest arc separates the others from the center of the circumcircle, whereas no arc on the left has this property.

\begin{proof}  We begin by stipulating some basic facts.  First, for any $J>0$ and $\theta_1,\hdots,\theta_{n-1}\in(0,2\pi)$ it is easy to see that any circle $C$ of radius $J$ contains a collection $\{x_0,\hdots,x_{n-1}\}$ such that the angle from $x_{i-1}$ to $x_i$ is $\theta_i$ for each $i>0$, and that this collection is cyclically ordered if and only if $\sum_{i=1}^{n-1} \theta_i < 2\pi$.  

Second, any such collection is determined up to isometry by the $\theta_i$: the orientation preserving isometry group of $\mathbb{H}^2$ acts transitively on circles of a fixed radius and includes all rotations around the center of $C$, and for another such collection $\{x_0',\hdots,x_{n-1}'\}$ on a circle $C'$ of radius $J$, the orientation-preserving isometry taking $C'$ to $C$ and $x_0'$ to $x_0$ takes $x_i'$ to $x_i$ for all $i$.

Finally, for any $x$ and $y$ at distance $d>0$ on a circle $C$ of radius $J$, $d\leq 2J$ by the triangle inequality, and the angle from $x$ to $y$ is $\theta(d,J)$ or $2\pi-\theta(d,J)$.  This follows upon applying Lemma \ref{angle function} to the triangle $T$ spanned by $x$, $y$ and the center $v$ of $C$.  The interior angle of $T$ at $v$ faces the counterclockwise arc from $x$ to $y$ in the former case, and the clockwise arc in the latter.

Now given $d_1,\hdots,d_{n-1}\in (0,\infty)$ and $J\geq \max\{d_i/2\}$, for any circle $C$ of radius $J$  the relation (\ref{LOS leader}) from the proof of Lemma \ref{angle function} implies that a collection $\{x_i\}\subset C$ with the property that the angle $\theta_i$ from $x_{i-1}$ to $x_i$ is $\theta(d_i,J)$ for each $i>0$ has $\mathit{dist}(x_{i-1},x_i) = d_i$ for each $i>0$.  This still holds if $\theta_{i_0} = 2\pi - \theta(d_{i_0},J)$ for some $i_0$, since in this case the interior angle of the triangle determined by $x_{i-1}$, $x_i$ and the center of $C$ faces the clockwise arc from $x_{i_0-1}$ to $x_{i_0}$ and hence still equals $\theta(d_{i_0},J)$.  Note also in this case that $\sum \theta_i < 2\pi$ if and only if $\theta(d_{i_0},J) > \sum_{i\neq i_0}\theta(d_i,J)$.

Now suppose $\{x_0,\hdots,x_{n-1}\}$ is a convex cyclic $n$-gon with side length collection $(d_0,\hdots,d_{n-1})$ on a hyperbolic circle $C$ of radius $J$, and for each $i>0$ let $\theta_i$ be the angle from $x_{i-1}$ to $x_i$.  From above we have $\theta_i = \theta(d_i,J)$ or $\theta_i = 2\pi - \theta(d_i,J)$ for each $i>0$, and $\sum\theta_i < 2\pi$.  Since $\theta(d_i,J)\in (0,\pi]$, $\theta_i$ can be $2\pi-\theta(d_i,J)$ for at most one $i>0$.  Moreover, if $i_0$ is such an index then $\sum_{i\neq i_0} \theta(d_i,J) < \theta(d_{i_0},J)$ as above.  In particular $d_{i_0}>d_i$ for each $i\neq i_0$ since one easily sees that $\theta(d,J)$ increases with $d$ for fixed $J$.

Below let us say $\{x_0,\hdots,x_{n-1}\}$ satisfies Case $n$ if $\theta_i<\pi$ for each $i$, and Case $s$ otherwise.  In Case $n$ the angle from $x_{n-1}$ to $x_0$ is $2\pi - \sum_{i=1}^{n-1} \theta(d_i,J)$, whereas in Case $s$ this angle is $\theta(d_{i_0},J) - \sum_{i\neq i_0} \theta(d_i,J)$.  Appealing to Lemma \ref{angle function} translates this angle measure into the functions $\ell^n$ and $\ell^s$ recorded above (in Cases $n$ and $s$, respectively).

Note that in case $n$, the triangle determined by $x_0$, $x_{n-1}$ and the circle center $v$ has interior angle $\theta = 2\pi - \sum_{i=1}^{n-1}\theta(d_i,J)$ at $v$ if $\theta\leq \pi$, but the interior angle here is $2\pi - \theta$ if $\theta > \pi$.  But this makes no difference in $\ell^n$ since $\sin x = \sin (\pi-x)$.\end{proof}

\begin{proposition}\label{omnibus2}  The functions $\ell^n$ and $\ell^s$ from Proposition \ref{omnibus1} are continuous on their domain 
$\left\{(J,d_1,\hdots,d_{n-1})\,|\, n\geq 3, d_i>0\ \mbox{for all}\ i,\ \mbox{and}\ J\geq \max\{d_i/2\}_{i=1}^{n-1}\right\},$ and smooth on its interior, with identical values on its frontier $\left\{(J,d_1,\hdots,d_{n-1})\,|\,J = \max\{d_i/2\}_{i=1}^{n-1}\right\}$.

Now fix $(d_1,\hdots,d_{n-1})\in(\mathbb{R}^+)^{n-1}$ for $n\geq 3$, let $D = \max\{d_i\}_{i=1}^{n-1}$ and, restricting $\ell^n$ and $\ell^s$ to $[D/2,\infty)\times\{(d_1,\hdots,d_{n-1})\}$, take them as functions of $J$.  Then $\frac{\partial \ell^n}{\partial J}(J)>0$ for all $J>J_0 = \min\{J\geq D/2\,|\,\sum_{i=1}^{n-1} \theta(d_i,J)\leq2\pi\}$, and $\lim_{J\to\infty} \ell^n(J) = 2\sinh^{-1}\left(\sum_{i=1}^{n-1} \sinh(d_i/2)\right)$.

For $i_0\in\{0,\hdots,n-1\}$ there exists $J\geq D/2$ such that $\theta(d_{i_0}, J)\geq \sum_{i\neq i_0} \theta(d_i,J)$ if and only if $d_{i_0} = D$ and $\sum_{i=1}^{n-1} \theta(d_i,D/2) \leq 2\pi$ (so in particular $J_0=D/2$).  If this is the case then $\{J\geq D/2\,|\,\theta(d_{i_0}, J)\geq \sum_{i\neq i_0} \theta(d_i,J)\}$ is an interval $I$ and $\frac{\partial\ell^s}{\partial J} (J) < 0$ for $J$ in its interior.

If $\sinh(d_{i_0}/2) \geq \sum_{i\neq i_0} \sinh(d_i/2)$ then $\sum_{i=1}^{n-1} \theta(d_i,D/2) \leq 2\pi$, $I$ as defined above is $[d_{i_0}/2,\infty)$, and $\lim_{J\to\infty} \ell^s(J) = 2\sinh^{-1}\left(\sinh(d_{i_0}/2)-\sum_{i\neq i_0}\sinh(d_i/2)\right)$.

If $\sinh(d_{i_0}/2) < \sum_{i\neq i_0} \sinh(d_i/2)$ but $\sum_{i=1}^{n-1} \theta(d_i,D/2) \leq 2\pi$ then $I = [d_{i_0}/2,J_1]$ for some $J_1<\infty$, and $\ell^s(J_1) = 0$.
\end{proposition}

\begin{remark}  Below, when a value of $(d_1,\hdots,d_{n-1})$ is clear from context we will frequently abbreviate $\ell^n(J,d_1,\hdots,d_{n-1})$ to $\ell^n(J)$ without comment; and likewise for $\ell^s$.\end{remark}

\begin{proof}  The continuity and smoothness properties of $\ell^n$ and $\ell^s$, and the description of their domain, follow directly from those of the functions involved in their definition.  In particular, see Lemma \ref{angle function} for those of $\theta(d,J)$.  That $\ell^n(J) = \ell^s(J)$ if $J = \max\{d_i/2\}$ follows from the fact that $\theta(d,d/2) = \pi$ and $\sin(\pi/2-x) = \sin(\pi/2+x)$ for any $x$.

We now fix $(d_1,\hdots,d_{n-1})\in(\mathbb{R}^+)^{n-1}$, let $D = \max\{d_i\}$ and, taking the restriction of $\ell^n$ to $[D/2,\infty)\times\{(d_1,\hdots,d_{n-1}\}$ as a function of $J$, prove its  properties described in the statement.  Taking $\theta_0(J) = \sum_{i=1}^{n-1}\theta(d_i,J)$, a derivative computation using Lemma \ref{angle function} gives:\begin{align}\label{d_0^n deriv}
  \frac{\partial}{\partial J} \ell^n(J) =  \frac{2\cosh J}{\cosh(\ell^n(J)/2)}\left[\sin(\theta_0(J)/2)-\cos(\theta_0(J)/2)\sum_{i=1}^{n-1} \tan(\theta(d_i,J)/2) \right] \end{align}

\begin{subclaim}\label{d_0^n increasing}For $J> J_0 = \min\{J\geq \max\{d_i/2\}\,|\,\sum_{i=1}^{n-1} \theta(d_i,J)\leq2\pi\}$, $\frac{\partial}{\partial J}\ell^n(J)>0$.\end{subclaim}

If $\theta_0(J_0)<\pi$ let $J_1 = J_0$; otherwise define $J_1>J_0$ by the equation $\theta_0(J_1) = \pi$.  If $J_1>J_0$ then for $J_0< J\leq J_1$ the derivative (\ref{d_0^n deriv}) is positive, since $\cos\left(\theta_0(J)/2\right)\leq 0$ for such $J$.  On $(J_1,\infty)$ we use:
\begin{claim}\label{add tans}  If $0< \sum_{i=1}^n \alpha_i <\pi/2$ then $\tan\left(\sum_{i=1}^n \alpha_i\right) > \sum_{i=1}^n \tan\left(\alpha_i\right)$\end{claim}
\begin{proof}[Proof of claim] For $0\leq\alpha\leq\beta$ with $\alpha+\beta<\pi/2$, an angle-addition identity gives:
$$ \tan(\alpha+\beta) = \frac{\tan\alpha+\tan\beta}{1-\tan\alpha\tan\beta} $$
Since $\alpha<\pi/2-\beta$ and sine is increasing on $(0,\pi/2)$ we have $\sin \alpha<\sin(\pi/2-\beta)=\cos\beta$.  Similarly $\cos\alpha>\sin\beta$, and so $0<\tan\alpha\tan\beta<1$.  Hence $\tan\alpha+\tan\beta < \tan(\alpha+\beta)$, and the claim follows by induction.\end{proof}

The claim implies that $\sum_{i=1}^{n-1} \tan(\theta(d_i,J)/2) < \tan(\theta_0(J)/2)$, hence that $\frac{\partial}{\partial J}\ell^n(J)>0$ on $(J_1,\infty)$.  This proves \ref{d_0^n increasing}.  We next show:

\begin{subclaim}\label{d_0^n limit} $\lim_{J\to\infty} \ell^n(J) = 2\sinh^{-1}\left(\sum_{i=1}^{n-1} \sinh(d_i/2)\right)$.\end{subclaim}

Again let $\theta_0(J)=\sum_{i=1}^{n-1} \theta(d_i,J)$.  The limit of $\sinh J \sin\left(\theta_0(J)/2\right)$ is of the form $\infty\cdot0$ by Lemma \ref{angle function}, so replacing $\sinh J$ by $\frac{1}{1/\sinh J}$ and applying l'H\^opital's rule gives:
$$ \lim_{J\to\infty} \sinh J \sin\left(\theta_0(J)/2\right) = 
  \lim_{J\to\infty} \cos\left(\theta_0(J)/2\right) \sum_{i=1}^{n-1} \frac{\sinh J\sinh(d/2)}{\sqrt{\sinh^2 J - \sinh^2(d/2)}} $$
Here we have used a formula from Lemma \ref{angle function}.  Evaluating the limit gives \ref{d_0^n limit}.

We now turn our attention to $\ell^s$.

\begin{claim}\label{ratio limit}  For positive real numbers $d_1,\hdots,d_{n-1}$ ($n\geq 3$) and any fixed $i_0$,
$$ \lim_{J\to\infty} \frac{\sum_{i\neq i_0} \theta(d_i,J)}{\theta(d_{i_0},J)} = \lim_{J\to\infty} \frac{\sum_{i\neq i_0} \frac{\partial}{\partial J} \theta(d_i,J)}{\frac{\partial}{\partial J} \theta(d_{i_0},J)} = \frac{\sum_{i\neq i_0} \sinh(d_i/2)}{\sinh(d_{i_0}/2)} $$
Moreover, if $d_{i_0}>d_i$ for all $i\neq i_0$ then $\left(\sum_{i\neq i_0} \frac{\partial}{\partial J} \theta(d_i,J)\right)/\frac{\partial}{\partial J} \theta(d_{i_0},J)$ is strictly increasing, and it is less than $1$ for every $J>d_{i_0}/2$ such that $\theta(d_{i_0},J) \geq \sum_{i\neq i_0}\theta(d_i,J)$.\end{claim}

\begin{proof}[Proof of claim]  Since each $\theta(d_i,J)\to 0$ as $J\to\infty$, the left-hand equation above follows from l'H\^opital's rule.  Directly substituting the formulas of Lemma \ref{angle function} gives:
$$ \frac{\sum_{i\neq i_0} \frac{\partial}{\partial J} \theta(d_i,J)}{\frac{\partial}{\partial J} \theta(d_{i_0},J)} =  \sum_{i\neq i_0} \frac{\sinh(d_i/2)}{\sinh(d_{i_0}/2)} \sqrt{\frac{\sinh^2 J - \sinh^2(d_{i_0}/2)}{\sinh^2 J - \sinh^2(d_i/2)} }$$
One thus easily verifies the limit computation.  Moreover, a slight rearrangment gives:
$$ \frac{\sum_{i\neq i_0} \frac{\partial}{\partial J} \theta(d_i,J)}{\frac{\partial}{\partial J} \theta(d_{i_0},J)} =  \sum_{i\neq i_0} \frac{\sinh(d_i/2)}{\sinh(d_{i_0}/2)} \sqrt{1 - \frac{\sinh^2 (d_{i_0}/2) - \sinh^2(d_{i}/2)}{\sinh^2 J - \sinh^2(d_i/2)}} $$
This makes clear that the derivative ratio is increasing.  To prove the final assertion we rewrite it using the tangent formula of Lemma \ref{angle function} and apply Claim \ref{add tans}.\end{proof}

\begin{subclaim}  If $\sinh(d_{i_0}/2) \geq\sum_{i\neq i_0} \sinh(d_i/2)$ then $\theta(d_{i_0},J) > \sum_{i\neq i_0} \theta(d_i,J)$ for all $J\geq d_{i_0}/2$, and $\lim_{J\to\infty} \ell^s(J) = 2\sinh^{-1}\left(\sinh(d_{i_0}/2)-\sum_{i\neq i_0}\sinh(d_i/2)\right)$.\end{subclaim}

If $\sinh(d_{i_0}/2) \geq\sum_{i\neq i_0} \sinh(d_i/2)$ then the derivative ratio $\left(\sum_{i\neq i_0} \frac{\partial}{\partial J} \theta(d_i,J)\right)/\frac{\partial}{\partial J} \theta(d_{i_0},J)$ is less than $1$ for all $J>d_{i_0}/2$, by Claim \ref{ratio limit}, so $\theta(d_{i_0},J)$ decreases more quickly than $\sum_{i\neq i_0}\theta(d_i,J)$ on the entire interval $(d_{i_0}/2,\infty)$.  It follows that if $\theta(d_{i_0},J)\leq\sum_{i\neq i_0}\theta(d_i,J)$ for any $J\geq d_{i_0}/2$ then for any $\epsilon>0$ there exists $\delta>0$ such that $\theta(d_{i_0},J)<\sum_{i\neq i_0}\theta(d_i,J)-\delta$ on $[J+\epsilon,\infty)$.  But this would violate the fact that both of these functions limit to $0$ (by Lemma \ref{angle function}), a contradiction.

The limit computation is an application of l'H\^opital's rule analogous to \ref{d_0^n limit}.

\begin{subclaim}\label{interval}  If $\theta(d_{i_0},J)\geq \sum_{i\neq i_0}\theta(d_i,J)$ for some $J\geq d_{i_0}/2$ then the set $I$ of all $J$ for which this inequality holds is an interval with left endpoint $d_{i_0}/2$.  If $\sinh(d_{i_0}/2) < \sum_{i\neq i_0} \sinh(d_i/2)$ then $I = [d_{i_0}2,J_1]$ is compact, with $\ell^s(J_1) = 0$.\end{subclaim}

These assertions both follow from Claim \ref{ratio limit}.  Its final assertion implies that $I$ is ``open to the left'' on $[d_{i_0}/2,\infty)$, so since it is also closed it is of the form $[d_{i_0}/2,\infty)$ or $[d_{i_0}/2, J_1]$ for some $J_1<\infty$.  But if $\sinh(d_{i_0}/2) < \sum_{i\neq i_0} \sinh(d_i/2)$ then the limit computation of Claim \ref{ratio limit} implies $\lim_{J\to\infty} \left(\sum_{i\neq i_0} \theta(d_i,J)\right)/ \theta(d_{i_0},J) > 1$, so $I$ must be of the latter form.  In this case $\ell^s(J_1) = 0$ by continuity, since $\ell^s(J)<0$ for all $J> J_1$.

It remains to show that $\partial\ell^s/\partial J<0$ in the interior of $I$.  Here let $\theta_0(J) = \theta(d_{i_0},J) - \sum_{i\neq i_0} \theta(d_i,J)$.  Then $\frac{\partial}{\partial J}\ell^s(J)$ has the following description:\begin{align*}
  \frac{2\cosh J}{\cosh(\ell^s(J)/2)}\left[ \sin(\theta_0(J)/2) 
     - \cos(\theta_0(J)/2)\left(\tan(\theta(d_{i_0},J)/2) - \sum_{i\neq i_0} \tan(\theta(d_i,J)/2)\right)\right]\end{align*}
One shows that $\tan(\theta(d_{i_0},J)/2) - \sum_{i\neq i_0} \tan(\theta(d_i,J)/2) > \tan(\theta_0/2)$ in two steps, first applying Claim \ref{add tans} to $\sum_{i\neq i_0} \tan(\theta(d_i,J)/2)$, then using the inequality $\tan(\alpha-\beta)<\tan\alpha-\tan\beta$ (which follows from the angle-addition identity for tangent).  This shows $\frac{\partial}{\partial J}\ell^s(J)<0$.\end{proof}

\begin{proposition}\label{existence uniqueness}  For positive real numbers $d_0,\hdots,d_{n-1}$ ($n\geq 3$) there is a cyclic $n$-gon with side length collection $(d_0,\hdots,d_{n-1})$ if and only if for each $i$, $\sinh(d_i/2) < \sum_{j\neq i} \sinh(d_j/2)$.  Two cyclic $n$-gons, with side length collections $(d_0,\hdots,d_{n-1})$ and $(d_0',\hdots,d_{n-1}')$, are isometric if and only if $d_i' = d_{\sigma(i)}$ for a cyclic permutation $\sigma$.  In particular, circumcircle radius is uniquely determined by $(d_0,\hdots,d_{n-1})$.  Moreover, it is symmetric as a function of $(d_0,\hdots,d_{n-1})$.\end{proposition}

\begin{proof}  We first prove existence.  For a fixed collection $d_0,\hdots,d_{n-1}$ ($n\geq 3$) of positive real numbers such that $\sinh(d_i/2) < \sum_{j\neq i} \sinh(d_j/2)$ for each $i$, we fix $i_0>0$ such that $d_{i_0}$ is maximal among $d_1,\hdots,d_{n-1}$ and consider three cases:\begin{enumerate}
  \item $\sum_{i=1}^{n-1} \theta(d_i,d_{i_0}/2) > 2\pi$
  \item $\sum_{i=1}^{n-1} \theta(d_i,d_{i_0}/2) \leq 2\pi$ but $\sinh(d_{i_0}/2) < \sum_{i=1,\ i\neq i_0}^{n-1} \sinh(d_i/2)$
  \item $\sinh(d_{i_0}/2) \geq \sum_{i=1,\ i\neq i_0}^{n-1}\sinh(d_i/2)$, hence $\sum_{i=1}^{n-1} \theta(d_i,d_{i_0}/2) \leq 2\pi$ by Proposition \ref{omnibus2}.\end{enumerate}

If $d_1,\hdots, d_{n-1}$ satisfies (1) then for $J_0$ as in Proposition \ref{omnibus2}, $\ell^n(J_0) = 0$.  Since $\ell^n$ limits to $2\sinh^{-1}\left(\sum_{i=1}^{n-1}\sinh(d_i/2)\right)$, and since $\sinh(d_0/2)<\sum_{i=1}^{n-1}\sinh(d_i/2)$, there exists $J> J_0$ such that $d_0 = \ell^n(J)$.  Since $J>J_0$, $\sum_{i=1}^{n-1}\theta(d_i,J) <2\pi$.  For this $J$, Proposition \ref{omnibus1} thus implies that a collection $x_0,\hdots,x_{n-1}$  arranged on a circle of radius $J$ such that the angle from $x_{i-1}$ to $x_i$ is $\theta(d_i,J)$ for each $i>0$ is a cyclic $n$-gon with side length collection $(d_0,\hdots,d_{n-1})$.

In case (2) let $D_0 = \ell^n(d_{i_0}/2) = \ell^s(d_{i_0}/2)$.  Proposition \ref{omnibus2} implies that the range of $\ell^s$ is $[0,D_0]$, and that of $\ell^n$ is $\left[D_0,2\sinh^{-1}\left(\sum_{i=1}^{n-1}\sinh(d_i/2)\right)\right)$, so as in case (1) there exists $J\geq d_{i_0}/2$ such that either $d_0 = \ell^n(J)$ or $\ell^s(J)$.  In the former instance we arrange $x_0,\hdots,x_{n-1}$ so that the angle from $x_i$ to $x_{i-1}$ is $\theta(d_i,J)$ for all $i$; and in the latter so that this holds for all $i$ but $i_0$, with the angle from $x_{{i_0}-1}$ to $x_{i_0}$ equal to $2\pi-\theta(d_{i_0},J)$.  Proposition \ref{omnibus1} again implies that $x_0,\hdots,x_{n-1}$ is a cyclic $n$-gon with the requisite side lengths.

Case (3) parallels case (2).  We must only note in addition that rearranging the inequality $\sinh(d_{i_0}/2) < \sum_{j\neq i_0} \sinh(d_j/2)$ (for $j$ between $0$ and $n-1$) gives $\sinh(d_0/2) > \sinh(d_{i_0}/2) - \sum_{i\neq i_0} \sinh(d_i/2)$ (for $i$ between $1$ and $n-1$).  This proves existence of a cyclic $n$-gon with side length collection $d_0,\hdots,d_{n-1}$ provided that $\sinh(d_i/2) < \sum_{j\neq i} \sinh(d_j/2)$ for each $i$.

If there is a cyclic $n$-gon with side length collection $d_0,\hdots,d_{n-1}$ then by Proposition \ref{omnibus1} either $d_0 = \ell^n(J)$ or $d_0 = \ell^s(J)$ for $\ell^n$ and $\ell^s$ as defined there, where $J$ is the circumcircle radius.  The range of $\ell^n$ has supremum $2\sinh^{-1}\left(\sum_{i=1}^{n-1} \sinh(d_i/2)\right)$, which is not attained, and if $\ell^s$ is defined it decreases from the minimum $D_0$, defined above, of $\ell^n$.  It follows that $\sinh(d_0/2) <\sum_{i=1}^{n-1} \sinh(d_i/2)$.

Moreover, if $d_{i_0}$ is maximal among $d_1,\hdots,d_{n-1}$ we claim that $\sinh(d_0/2) > \sinh(d_{i_0}/2) - \sum_{i\neq i_0} \sinh(d_{i_0}/2)$.  If $d_1,\hdots,d_{n-1}$ satisfies case (3) above this holds because the range of $\ell^s$ is bounded below by $2\sinh^{-1}\left(\sinh(d_{i_0}/2) - \sum_{i\neq i_0} \sinh(d_i/2)\right)$ by Proposition \ref{omnibus2}, and that of $\ell^n$ by $D_0$.  In case (2) the right-hand side of the claimed inequality is less than $0$ by hypothesis.  This also holds in case (1), since by Proposition \ref{omnibus2} $\sum_{i=1}^{n-1} \theta(d_i,d_{i_0}/2) \leq 2\pi$ if $\sinh(d_{i_0}/2) \geq \sum_{i\neq i_0}\sinh(d_i/2)$.  The claim is thus proved, and it follows that $\sinh(d_i/2) <\sum_{j\neq i} \sinh(d_j/2)$ for all $i$.

For a given $d_0,\hdots,d_{n-1}$ with $\sinh(d_i/2)\leq \sum_{j\neq i} \sinh(d_j/2)$ for all $i$, we claim there is only one $J$ such that $d_0 = \ell^n(J)$ or $\ell^s(J)$; i.e. that circumcircle radius is uniquely determined by the side length collection.  To this end we note first that cases (1), (2) and (3) above are mutually exclusive.  Moreover, $\ell^n$ and $\ell^s$ are each strictly monotone, and their ranges intersect only at the single point $D_0 = \ell^n(d_{i_0}/2) = \ell^s(d_{i_0}/2)$.  The claim follows.

It is further evident by inspecting the definitions of $\ell^n$ and $\ell^s$ that they are symmetric in $d_1,\hdots,d_{n-1}$, so the circumcircle radius $J$ determined by $d_0,\hdots,d_{n-1}$ is invariant under any reordering of the final $n-1$ side lengths.  We claim it is also invariant under cyclic permutations of $(d_0,\hdots,d_{n-1})$. 

The key here is that Proposition \ref{omnibus2} also applies to the restrictions of $\ell^n$ and $\ell^s$ to $[D_1/2,\infty)\times\{(d_2,\hdots,d_{n-1},d_0)\}$, where $D_1 = \max\{d_i\}_{i\neq 1}$.  The arguments above thus imply that $d_1 = \ell^n(J_1,d_2,\hdots,d_{n-1},d_0)$ or $d_1 = \ell^s(J_1,d_2,\hdots,d_{n-1},d_0)$ for a unique $J_1\in [D_1/2,\infty)$.  That $J_1 = J$ follows from uniqueness, with a little case-matching.

If $d_0 = \ell^n(J,d_1,\hdots,d_{n-1})$ then plugging into the definition of $\theta(d,J)$ from Lemma \ref{angle function} gives:
$$ \theta(d_0,J) = \sum_{i=1}^{n-1} \theta(d_i,J) \quad\mbox{or}\quad 
  \theta(d_0,J) = 2\pi - \sum_{i=1}^{n-1} \theta(d_i,J) $$
The ambiguity arises because the inverse sine takes values in $[-\pi/2,\pi/2]$: for $x\in[0,\pi]$, $\sin^{-1}(\sin x) = x$ if $x\leq \pi/2$ and $\pi-x$ otherwise.  The former case thus occurs if and only if $\sum_{i=1}^{n-1} \theta(d_i,J) \leq \pi$; the latter if and only if this sum lies between $\pi$ and $2\pi$.  Rearranging the former case gives $\theta(d_1,J) = \theta(d_0,J) - \sum_{i=2}^{n-1} \theta(d_i,J) \in (0,\pi)$, so $d_1 = \ell^s(J,d_2,\hdots,d_{n-1},d_0)$.  In the latter we have $\theta(d_1,J) = 2\pi - \sum_{i\neq 1} \theta(d_i,J)$, so $d_1 = \ell^n(J,d_2,\hdots,d_{n-1},d_0)$.  (Here we use that $\sin (\pi-x) = \sin x$.)  In either case we conclude from uniqueness that $J=J_1$.

If $d_0 = \ell^s(J,d_1,\hdots,d_{n-1})$ then $\theta(d_0,J) = \theta(d_{i_0},J) - \sum_{i>0,i\neq i_0}\theta(d_i,J)$, where $d_{i_0}$ is maximal among $d_1,\hdots,d_{n-1}$.  There are two cases, handled similarly to the above: if $i_0 = 1$ then $d_1 = \ell^n(J,d_2,\hdots,d_{n-1},d_0)$, otherwise $d_1 = \ell^s(J,d_2,\hdots,d_{n-1},d_0)$, and in each case we again conclude $J=J_1$.  The claim thus follows, and also that $J = J(d_0,\hdots,d_{n-1})$ is symmetric.

We finally show that cyclic $n$-gons with side length collections $(d_0,\hdots,d_{n-1})$ and $(d_0',\hdots,d_{n-1}')$ are isometric if and only if $d_i' = d_{\sigma(i)}$ for each $i$, where $\sigma$ is a fixed cyclic permutation.  The ``only if'' direction is in Remark \ref{important}, so let us suppose for cyclic  $n$-gons $\{x_i\}$ and $\{x_i'\}$ with respective side length collections above that there is a cyclic permutation $\sigma$ so that $d_i' = d_{\sigma(i)}$ for all $i$.  Upon cyclically reordering the $x_i$ we may assume that $d_i' = d_i$ for each $i$ (again see Remark \ref{important}).  We have already showed that $\{x_i\}$ and $\{x_i'\}$ have the same circumcircle radius $J$, and moreover that $d_0$ is only one of $\ell^n(J)$ or $\ell^s(J)$.  Proposition \ref{omnibus1} now implies that $\{x_i\}$ and $\{x_i'\}$ are isometric.\end{proof}

Having established existence and uniqueness of cyclic polygons by carefully analyzing the functions $\ell^n$ and $\ell^s$, we will use the following lemma to avoid having to further consider $\ell^s$.

\begin{lemma}\label{reduce}  Suppose a cyclic $n$-gon with side length collection $(d_0,\hdots,d_{n-1})$, $n\geq 3$, and circumcircle radius $J$ has $d_0 = \max\{d_i\}_{i=0}^{n-1}$.  Then the following hold: $d_0 = \ell^n(J,d_1,\hdots,d_{n-1})$, for $\ell^n$ as defined in Proposition \ref{omnibus1}; and $J > D/2$, for $D = \max\{d_i\}_{i=1}^{n-1}$.  Furthermore, both assertions hold for all side length collections of cyclic $n$-gons near $(d_0,\hdots,d_{n-1})$ in $\mathbb{R}^n$.\end{lemma}

\begin{proof}By Proposition \ref{omnibus1}, either $d_0 = \ell^n(J,d_1,\hdots,d_{n-1})$ or $d_0 = \ell^s(J,d_1,\hdots,d_{n-1})$ as defined there.  But $\ell^s$ is decreasing in $J$ and defined only for $J\geq D/2$ by Proposition \ref{omnibus2}, where $D= \max\{d_i\}_{i=1}^{n-1} \leq d_0$ by hypothesis.  We calculate:
$$ \ell^s(D/2,d_1,\hdots,d_{n-1}) < 2\sinh^{-1}\left[\sinh(D/2)\sin\left(\frac{1}{2}\theta(D,D/2)\right)\right] = D $$
Since $\ell^s$ is continuous this inequality holds on a neighborhood of $(d_1,\hdots,d_{n-1})$ in $\mathbb{R}^n$, and it follows that there is a neighborhood of $(d_0,\hdots,d_{n-1})$ such that $d_0 = \ell^n(J)$ for each side length collection in this neighborhood.

By Proposition \ref{omnibus2} we have $\ell^n(D/2,d_1,\hdots,d_{n-1}) = \ell^s(D/2,d_1,\hdots,d_{n-1}) < D$, so since $\ell^n$ is continuous this inequality also holds on an entire neighborhood of $(d_0,\hdots,d_{n-1})$.  It implies that $J>D/2$ here.\end{proof}

\subsection{Parametrizing cyclic polygons}

Having proved the existence and uniqueness of cyclic hyperbolic $n$-gons, in this section we will first parametrize them by side length using an open set in $(\mathbb{R}^+)^n$.  Then we will describe some subspaces of this set corresponding to subclasses of cyclic polygons that turn out to be geometrically significant.

\begin{corollary}\label{ACn}  For $n\geq 3$, the collection\begin{align*}
   \mathcal{AC}_n & = \left\{(d_0,\hdots,d_{n-1}) \in (\mathbb{R}^+)^n\,|\, \sinh (d_i/2) < \sum_{j\neq i} \sinh(d_j/2) \ \mbox{for each}\ i \in \{0,\hdots,n-1\}\right\} \end{align*}
parametrizes isometry classes of marked cyclic $n$-gons in $\mathbb{H}^2$ by their side length collections, where a \mbox{\rm marking} is a choice of vertex to label $x_0$.\end{corollary}

This is a direct consequence of Proposition \ref{existence uniqueness}, since the marking takes care of the cyclic ambiguity of side length collections.

\begin{proposition}\label{Cn and BCn}  For $n\geq 3$, a cyclic $n$-gon $\{x_0,\hdots,x_{n-1}\}$ is \mbox{\rm centered} if the angle from $x_{i-1}$ to $x_i$ is less than $\pi$ for each $i$; or equivalently, if $\sum_{i\neq i_0} \theta(d_i,J) >\pi$ where $(d_0,\hdots,d_{n-1})$ is its side length collection, with $d_{i_0}$ maximal, and $J$ is its circumcircle radius. The collection\begin{align*}
   \mathcal{C}_n & = \left\{(d_0,\hdots,d_{n-1}) \in (\mathbb{R}^+)^n\,|\,\sum_{i=0}^{n-1}\theta(d_i,d_{i_0}/2) > 2\pi,\ \mbox{where}\ d_{i_0}=\max\{d_i\}_{i=0}^{n-1} \right\} \subset \calAC_n\end{align*}
parametrizes marked, centered $n$-gons in $\mathbb{H}^2$ up to hyperbolic isometry.  The collection\begin{align*}
  \mathcal{BC}_n & = \left\{(d_0,\hdots,d_{n-1})\in (\mathbb{R}^+)^n\,|\,\sum_{i=0}^{n-1}\theta(d_i,d_{i_0}/2) = 2\pi,\ \mbox{where}\ d_{i_0}=\max\{d_i\}_{i=0}^{n-1} \right\} \subset \calAC_n \end{align*}
parametrizes marked, semicyclic $n$-gons in $\mathbb{H}^2$ up to hyperbolic isometry.  Here a cyclic $n$-gon is \mbox{\rm semicyclic} if its circumcircle radius $J$ is $d_{i_0}/2$; or equivalently if $\sum_{i\neq i_0}\theta(d_i,J) = \pi$.\end{proposition}


\begin{proof}  We first comment on the equivalent definitions of centeredness.  Recall from Proposition \ref{omnibus1} that for each $i$ that the angle from $x_i$ to $x_{i-1}$ is either $\theta(d_i,J)$, as defined in Lemma \ref{angle function}, or $2\pi - \theta(d_i,J)$.  For each $i$ such that this angle is less than $\pi$ it must equal $\theta(d_i,J)$ (recall that $\theta(d,J) \in [0,\pi]$ for all possible $d$ and $J$).  If this holds for all $i$ then since these angles sum to $2\pi$ it follows in particular that $\sum_{i\neq i_0} \theta(d_i,J) = 2\pi - \theta(d_{i_0},J) >\pi$.  

On the other hand if $\sum_{i\neq i_0} \theta(d_i,J) > \pi$ for $i_0$ such that $d_{i_0}$ is maximal then the same inequality holds for any other $i_0$, since $\theta(d,J)$ increases in $d$ for fixed $J$.  It follows that for each $i$, $(2\pi - \theta(d_{i},J))+\sum_{j\neq i}\theta(d_j,J) > 2\pi$, so the angle from $x_{i-1}$ to $x_i$ is $\theta(d_i,J) <\pi$. 

Now suppose a centered or semicyclic $n$-gon has side length collection $(d_0,\hdots,d_{n-1})$ and circumcircle radius $J$.  We may assume $d_0$ is maximal among the $d_i$, since circumcircle radius is symmetric in side lengths and the centeredness criterion is invariant under their permutation.   Then by Lemma \ref{reduce}, $d_0 = \ell^n(J,d_1,\hdots,d_{n-1})$ so manipulating the definition of $\ell^n$ gives:
$$ \sinh(d_0/2) = \sinh J \sin\left(\frac{1}{2}\sum_{i=1}^{n-1} \theta(d_i,J)\right) \leq \sinh J $$
It follows that $J \geq d_0/2$ (with equality here and above if and only if the $n$-gon is semicyclic), hence since $\theta(d,J)$ decreases in $J$ that: 
$$\sum_{i=0}^{n-1} \theta(d_i,d_{i_0}/2) = \pi+\sum_{i=1}^{n-1} \theta(d_i,d_{i_0}/2) \geq \pi+\sum_{i=1}^{n-1} \theta(d_i,J) \geq 2\pi$$
Again, equality holds if and only if the original $n$-gon is semicyclic, and we find that the side length collection of a centered or semicyclic collection lies in $\calc_n$ or $\calBC_n$, respectively.

Now suppose on the other hand that $(d_0,\hdots,d_{n-1})\in\calc_n\cup\calBC_n$ has $d_0$ maximal.  Then since $\sum_{i=0}^{n-1} \theta(d_i,d_0/2) \geq 2\pi$ there exists a unique $J\geq d_0/2$ such that $\sum_{i=0}^{n-1} \theta(d_i,J) = 2\pi$, with $J = d_0/2$ if and only if $\sum_{i=0}^{n-1} \theta(d_i,d_0/2) = 2\pi$.

For this $J$, directly substituting $2\pi - \theta(d_0,J)$ for $\sum_{i=1}^{n-1} \theta(d_i,J)$ in the definition of $\ell^n$ and applying $\theta$ we find that $\theta(d_0,J) = \theta(\ell^n(J),J)$, so $d_0 = \ell^n(J)$.  Therefore $J$ is the circumcircle radius of a cyclic $n$-gon with side length collection $(d_0,\hdots,d_{n-1})$.  If $J = d_0/2$ then the $n$-gon in question is semi-cyclic; otherwise $\theta(d_0,J) < \pi$ so $\sum_{i=1}^{n-1} \theta(d_i,J) > \pi$ and it is centered.

We note further that for $(d_0,\hdots,d_{n-1})\in\calc_n\cup\calBC_n$, since $d_0 = \ell^n(J)$ and $\ell^n$ increases toward an asymptote of $2\sinh^{-1}\left(\sum_{i=1}^{n-1}\sinh(d_i/2)\right)$, that $(d_0,\hdots,d_{n-1})\in\calAC_n$. 

We finally comment on equivalence of the definitions of semicyclicity.  For a cyclic $n$-gon with side length collection $(d_0,\hdots,d_{n-1})$ such that $d_0$ is maximal, so that $d_0 = \ell^n(J)$ by Lemma \ref{angle function}, directly substituting either $J = d_{i_0}/2$ or $\sum_{i=1}^{n-1} \theta(d_i,J) = \pi$ into the formula for $\ell^n(J)$ allows one to conclude the other condition.\end{proof}

\begin{proposition}\label{submanifold}  For each $n\geq 3$, $\calBC_n$ is the frontier of $\calc_n$ in $\calAC_n$.  It is the orbit of $\mathrm{graph}(b_0)\doteq \{(b_0(\bd),\bd)\,|\,\bd\in(\mathbb{R}^+)^{n-1}\}$ under cyclic  permutation of entries, where $b_0\co(\mathbb{R}^+)^{n-1}\to \mathbb{R}$ is defined by:
$$ \sum_{i=1}^{n-1} \theta(d_i,b_0(\bd)/2) = \pi$$
This function is symmetric, smooth, and strictly increasing in each variable.  It further satisfies $b_0(d_1,\hdots,d_{n-1}) > \max\{d_i\}_{i=1}^{n-1}$.  In particular, $\calBC_n$ has $n$ connected components.\end{proposition}

\begin{proof}  Proposition \ref{Cn and BCn} directly implies that $\calBC_n$ is the frontier of $\calc_n$ in $\calAC_n$, since the quantities involved in the definitions of these sets vary continuously over $\calAC_n$.  

If $(d_0,\hdots,d_{n-1})\in\calBC_n$ has maximal entry $d_{i_0}$ then $\sum_{i\neq i_0} \theta(d_i,d_{i_0}/2) = \pi$, since $\theta(d_{i_0},d_{i_0}/2) = \pi$.  For $D = \max\{d_i\}_{i\neq i_0}$, $\sum_{i\neq i_0} \theta(d_i,D/2) > \pi$, so since $\theta(d,J)$ strictly decreases in $J$ for each $(d_0,\hdots,\widehat{d_{i_0}},\hdots,d_{n-1})\in(\mathbb{R}^+)^{n-1}$, the equation $\sum_{i\neq i_0}\theta(d_i,d_{i_0}/2) = \pi$ uniquely determines $d_{i_0} > D$.

The above implies in particular that $b_0\co(\mathbb{R}^+)^{n-1}\to\mathbb{R}$ is uniquely determined by its defining equation, that $b_0(d_1,\hdots,d_{n-1}) > \max\{d_i\}$, and that $\calBC_n$ is the orbit of $\mathrm{graph}(b_0)$ under cyclic permutation of entries.  Furthermore since $b_0 > \max\{d_i\}$ no distinct translates of $\mathrm{graph}(b_0)$ under this action intersect each other.  We will thus finish the proof by showing that $b_0$ is smooth and strictly increasing in each variable.

This is a direct application of the implicit function theorem, which shows that for each $i$,
$$ \frac{\partial}{\partial d_i} b_0(d_1,\hdots,d_{n-1}) = -\frac{\frac{\partial\theta}{\partial d}(d_i,b_0(d_1,\hdots,d_{n-1}))}{\sum_{i=1}^{n-1} \frac{\partial\theta}{\partial J}(d_i,b_0(d_1,\hdots,d_{n-1}))} $$
Since $\theta$ strictly increases in $d$ and decreases in $J$, the result follows.
\end{proof}

\subsection{Smoothness}  Having parametrized cyclic $n$-gons by the space $\calAC_n$, we turn our attention to describing associated geometric quantities as functions on $\calAC_n$: circumcircle radius in Proposition \ref{smooth J}, diagonal lengths in Corollary \ref{diagonals}, and what turns out to be area in Corollary \ref{area function}.  The main result in each case is that the function in question is smooth, though Proposition \ref{smooth J} also describes some qualitative features of the circumcircle radius.

\begin{lemma}\label{J derivs} For $n\geq 3$, suppose $\bd = (d_0,\hdots,d_{n-1})\in\calAC_n$ has $d_0\geq d_i$ for all $i$.  The function $J\co\calAC_n\to\mathbb{R}^+$ that records circumcircle radius is smooth at $\bd$, and:\begin{align*}
  \frac{\partial J}{\partial d_i}(\bd) = \frac{1}{2}\frac{\pm \cosh(d_i/2)\tanh J \sqrt{\frac{\sinh^2 J - \sinh^2(d_0/2)}{\sinh^2 J - \sinh^2(d_i/2)}}}{\sinh(d_0/2) \pm \sum_{j=1}^{n-1}\sinh(d_j/2)\sqrt{\frac{\sinh^2 J - \sinh^2(d_0/2)}{\sinh^2 J - \sinh^2(d_j/2)}}}  \end{align*}
Above, each ``$\pm$'' should be taken as ``$+$'' if $\bd\in\calc_n$ and ``$-$'' otherwise, with one exception: if $\bd\in\calAC_n-\calc_n$ and $i=0$, the numerator is $\cosh(d_0/2)\tanh J$.\end{lemma}

\begin{proof}  Since $d_0$ is maximal, Lemma \ref{reduce} gives a neighborhood $U$ of $\bd$ in $\mathcal{AC}_n$ on which the equation $d_0 = \ell^n(J,d_0,\hdots,d_{n-1})$ holds everywhere, for $J = J(d_0,\hdots,d_{n-1})$.  We will thus obtain the desired conclusions by applying the implicit function theorem to $F(d_0,\hdots,d_{n-1},J) = \ell^n(J,d_1,\hdots,d_{n-1}) - d_0$.  Since $J > D/2$ by Lemma \ref{reduce}, where $D = \max\{d_i\}_{i=1}^{n-1}$, using Lemma \ref{angle function} and inspecting the definition of $\ell^{n}$ we find that $F$ is smooth near $(\bd,J)$.    Its derivative vector is $\left(-1,\frac{\partial}{\partial d_1}\ell^n,\hdots,\frac{\partial}{\partial d_{n-1}}\ell^n, \frac{\partial}{\partial J} \ell^n \right)$.

Since $J(\bd) > J_0$ as defined in Proposition \ref{omnibus2}, that result gives $\frac{\partial}{\partial J} \ell^n(J(\bd))>0$.  The implicit function theorem thus guarantees that $J$ is smooth at $\bd$, and the following relations hold:\begin{align}
  \label{zero} 0 &= \frac{\partial J}{\partial d_0}\frac{\partial\ell^n}{\partial J} - 1 \\
  \label{bigger} 0 &= \frac{\partial J}{\partial d_i}\frac{\partial\ell^n}{\partial J} + \frac{\partial \ell^n}{\partial d_i} \qquad (i>0)
\end{align}
We will substitute for $\frac{\partial}{\partial J} \ell^n$ in the relations above using (\ref{d_0^n deriv}).  In that formula, $\theta_0(J) = \sum_{i=1}^{n-1} \theta(d_i,J)$, and using the definition of $\ell^n$ and the fact that $d_0 = \ell^n(J(\bd))$ we obtain:\begin{align*}
  & \sin(\theta_0(J)/2) = \frac{\sinh(d_0/2)}{\sinh J} &
  & \cos(\theta_0(J)/2) = \pm \frac{\sqrt{\sinh^2 J - \sinh^2(d_0/2)}}{\sinh J} \end{align*}
The ``$\pm$'' above is a ``$-$'' if $\theta_0(J) >\pi$ and a ``$+$'' otherwise; by Proposition \ref{Cn and BCn} these cases respectively correspond to $\bd\in\calc_n$ and $\bd\in\calAC_n-\calc_n$.  We further substitute $\sinh(d_i/2)/\sqrt{\sinh^2 J - \sinh^2(d_i/2)}$ for $\tan(\theta(d_i,J)/2)$, yielding:
$$ \frac{\partial}{\partial J}\ell^n = \frac{2\coth J}{\cosh(d_0/2)}\left[\sinh(d_0/2) \pm \sum_{i=1}^{n-1} \sinh(d_i/2)\sqrt{\frac{\sinh^2 J - \sinh^2(d_0/2)}{\sinh^2 J - \sinh^2(d_i/2)}}\right] $$
Here the ``$\pm$'' is a ``$+$'' if $\bd\in\calc_n$ and ``$-$'' otherwise.  With the same dichotomy below,  a similar calculation gives:
$$ -\frac{\partial \ell^n}{\partial d_i} = \pm \frac{\cosh(d_i/2)}{\cosh(d_0/2)}\sqrt{\frac{\sinh^2 J - \sinh^2(d_0/2)}{\sinh^2 J-\sinh^2(d_i/2)}} $$
The formula from the statement thus follows from (\ref{zero}) and (\ref{bigger}).\end{proof}

The bounds below distill a useful consequence of the formulas of Lemma \ref{J derivs}.

\begin{proposition}\label{smooth J}\SmoothJ\end{proposition}

\begin{proof}  $J$ is symmetric by Proposition \ref{existence uniqueness}, so Lemma \ref{J derivs} actually implies that it is smooth on all of $\calAC_n$.  To obtain the bounds above we will fix $\bd = (d_0,\hdots,d_{n-1})\in\mathcal{AC}_n$ and assume that $d_0$ is maximal among the $d_i$ by rearranging if necessary.  Then we will use the equations (\ref{zero}) and (\ref{bigger}) from the proof of Lemma \ref{J derivs}, which hold with this hypothesis.

Applying (\ref{d_0^n deriv}) from the proof of Proposition \ref{omnibus2} to (\ref{zero}), with $\theta_0(J) = \sum_{i=1}^{n-1} \theta(d_i,J)$, yields:
$$ \frac{\partial J}{\partial d_0} = \frac{1}{\frac{\partial}{\partial J}\ell^n(J)}
  = \frac{\cosh(d_0/2)}{2\cosh J}\ \frac{1}{\sin(\theta_0(J)/2)-\cos(\theta_0(J)/2)\sum_{i=1}^{n-1} \tan(\theta(d_i,J)/2)} $$
Note that if $\theta_0(J) = \pi$ then $d_0 = 2J$ and $\partial J/\partial d_0 = 1/2$.  We accordingly divide the remaining cases into two subcases: $\theta_0(J) < \pi$ and $\theta_0(J) > \pi$.

If $\theta_0(J)<\pi$ then since $\cos(\theta_0(J)/2)>0$ the denominator of the right-hand quantity above is less than $\sin(\theta_0(J)/2)$.  By the definition of $\ell^n$, $\sin (\theta_0(J)/2) = \sinh(\ell^n(J)/2)/\sinh(J)$, so since $\ell^n(J) = d_0$ we have the following inequality:
$$ \frac{\partial J}{\partial d_0} > \frac{\cosh(d_0/2)}{2\cosh J} \frac{\sinh J}{\sinh(d_0/2)} = \frac{\tanh(J)}{2\tanh (d_0/2)} $$
Since the hyperbolic tangent is increasing and $J > d_0/2$ it follows that $\partial J/\partial d_0 > 1/2$.  This is the case $\bd\in\mathcal{AC}_n - (\mathcal{C}_n \cup\mathcal{BC}_n)$ by Proposition \ref{Cn and BCn}.

If $\theta_0(J) > \pi$ then since $\cos(\theta_0(J)/2) < 0$ and $\tan x > x$ on $(0,\pi/2)$ we have 
$$\sin(\theta_0(J)/2)-\cos(\theta_0(J)/2)\sum_{i=1}^{n-1} \tan(\theta(d_i,J)/2) > \sin(\theta_0(J)/2)-\cos(\theta_0(J)/2)\cdot \pi/2 $$
Since $\sin x - \cos x > 1$ on $(\pi/2,\pi)$ we have $\partial J/\partial d_0 < \cosh(d_0/2)/(2\cosh J) < 1/2$ in this case, which by Proposition \ref{Cn and BCn} is $\bd\in\mathcal{C}_n$.

For $i>1$, solving (\ref{bigger}) for $\partial J/\partial d_i$ gives $\partial J/\partial d_i = -\frac{\frac{\partial}{\partial d_i} \ell^n(J)}{\frac{\partial}{\partial J} \ell^n (J)}$.  Let us record the numerator of this quantity:
$$ \frac{\partial}{\partial d_i} \ell^n(J) = \frac{\sinh J\cos\left(\theta_0(J)/2\right)}{\cosh(\ell^n(J)/2)}\ \frac{\partial}{\partial d} \theta(d_i,J) = \frac{\cosh(d_i/2)}{\cosh(d_0/2)}\, \frac{\cos(\theta_0(J)/2)}{\cos(\theta(d_i,J)/2)} $$
The second equation above follows from the definition in Lemma \ref{angle function} after taking a partial derivative with respect to $d$.  All terms here are positive except possibly $\cos(\theta_0(J)/2)$, so $\partial J/\partial d_i$ is negative if $\theta_0(J) < \pi$ (i.e.~$\bd\in\mathcal{AC}_n - (\mathcal{C}_n \cup\mathcal{BC}_n)$) and positive if $\theta_0(J) > \pi$ ($\bd\in\calc_n$).  

For fixed $J$, $\theta(d,J)$ increases in $d$ on $(0,2J)$, so $f(d) = \cosh(d/2)/\cos(\theta(d,J)/2)$ is also increasing.  Thus for $i$ and $j$ at least one we have $|\partial J/\partial d_i| > |\partial J/\partial d_j|$ if $d_i > d_j$.  We claim that also $|\cosh(d_0/2)/\cos(\theta_0(J)/2)| \geq f(d_i)$.  Manipulating the definition of $\ell^n(J)$ gives:
$$ \cos(\theta_0(J)/2) = \pm\frac{\sqrt{\sinh^2 J - \sinh^2(d_0/2)}}{\sinh(J)} $$
By Lemma \ref{angle function} $\cos(\theta(d_i,J)/2)$ satisfies the same equation on taking an absolute value and substituting $d_i$ for $d_0$ on the right-hand side.  Since $d_0\geq d_i$ we have $|\cos(\theta_0(J)/2)| < \cos(\theta(d_i,J)/2)$, and the claim follows.  This implies the result.
\end{proof}

\begin{corollary}\label{diagonals}  For any $n\geq 4$ and $i,j\in\{0,\hdots,n-1\}$ such that $|i-j|\geq 2\  (\mathrm{mod}\ 4)$, defining $\ell_{i,j}(d_0,\hdots,d_{n-1}) = \mathit{dist}(x_i,x_j)$, where $\{x_0,\hdots,x_{n-1}\}$ is a cyclic $n$-gon with side length collection $(d_0,\hdots,d_{n-1})$, yields a smooth function $\calAC_n\to\mathbb{R}^+$.\end{corollary}

\begin{proof}  Let us suppose without loss of generality that $i < j$, fix $(d_0,\hdots,d_{n-1})\in\calAC_n$, and let $\{x_0,\hdots,x_{n-1}\}$ be a cyclic $n$-gon with side length collection $(d_0,\hdots,d_{n-1})$ and circumcircle radius $J = J(d_0,\hdots,d_{n-1})$.  For $i_0$ such that $d_{i_0}$ is maximal, let us also assume for now that $i_0\notin \{i+1,\hdots,j\}$.  Then by Proposition \ref{omnibus1} the angle from $x_{k-1}$ to $x_k$ is $\theta(d_k,J)<\pi$ for each $k \in \{i+1,\hdots,j\}$.

Applying Proposition \ref{omnibus1}, mutatis mutandis, to the cyclically ordered collection $\{x_i,\hdots,x_j\}$ now gives:
$$\ell_{i,j}(d_0,\hdots,d_{n-1}) = 2\sinh^{-1}\left[\sinh J\sin\left(\frac{1}{2}\sum_{k=i+1}^j \theta(d_k,J)\right)\right]$$
We claim that the same formula holds in a neighborhood of $(d_0,\hdots,d_{n-1})$, whence smoothness of $\ell_{i,j}$ follows from smoothness of $J$ and of the functions above.  If it did not hold then for some $i_1\in\{i+1,\hdots,j\}$ there would be points $(d_0',\hdots,d_{n-1}')$ arbitrarily close to $(d_0,\hdots,d_{n-1})$ such that for a cyclic $n$-gon $\{x_0',\hdots,x_{n-1}'\}$ on a circle of radius $J'$ with side length collection $(d_0',\hdots,d_{n-1}')$, the angle from $x_{i_1-1}$ to $x_i$ was $2\pi-\theta(d_{i_1}',J')$.

If this were the case then we would have $\theta(d_{i_1}',J') = \sum_{i\neq i_1} \theta(d_i',J')$ for any such $(d_0',\hdots,d_{n-1}')$, which upon taking a limit would also hold for $(d_0,\hdots,d_{n-1})$.  But this would contradict the fact that $d_{i_0}$ is maximal among the $d_i$.  This proves the claim and hence the establishes the result in the case $i_0\notin\{i+1,\hdots,j\}$.

The case $i_0\in\{i+1,\hdots,j\}$ is analogous to the other, but in the formula for $\ell_{i,j}$ above we replace the sum over $\{i+1,\hdots,j\}$ with one over $\{j+1,\hdots,n-1,0,\hdots,i\}$.\end{proof}

\begin{lemma}\label{triangle area}  Let $U$ be the set of points in $(\mathbb{R}^+)^3$ such that each coordinate is less than the sum of the others.  For $(a,b,c)\in U$, define:
$$ \alpha(a,b,c) = \cos^{-1}\left(\frac{\cosh b\cosh c - \cosh a}{\sinh b\sinh c}\right) $$
The function $A(a,b,c) = \pi - \alpha(a,b,c) - \alpha(b,c,a) - \alpha(c,a,b)$ is smooth on $U$.  It records the hyperbolic area of a triangle in $\mathbb{H}^2$ with sides of length $a$, $b$ and $c$.\end{lemma}

\begin{proof}  Since the inverse cosine is smooth on $(-1,1)$, the only thing to note is that the quantity in parentheses above lies in this interval; i.e.~that:
$$ -\sinh b \sinh c < \cosh b \cosh c - \cosh a < \sinh b \sinh c $$
This follows from the ``angle sum'' identity for hyperbolic cosine and the fact that $(a,b,c)\in U$.  The hyperbolic law of cosines (see \cite[Th.~3.5.3]{Ratcliffe}) implies that $\alpha(a,b,c)$ is the angle of a hyperbolic triangle with sides of length $a$, $b$ and $c$ at the vertex opposite the side with length $a$.  That $A(a,b,c)$ measures the area of such a triangle follows from the Gauss--Bonnet formula for hyperbolic area, see Theorem 3.5.5 of \cite{Ratcliffe}.\end{proof}

\begin{corollary}\label{area function}  For $\bd= (d_0,d_1,d_2)\in\calAC_3$ define $D_0(\bd) = A(d_0,d_1,d_2)$, where $A$ is as defined in Lemma \ref{triangle area}.  For $\bd = (d_0,d_1,d_2,d_3)\in\calAC_4$ define: 
$$D_0(\bd) = A(d_0,d_1,\ell_{0,2}(\bd)) + A(\ell_{0,2}(\bd),d_2,d_3)$$
Here $\ell_{i,j}$ is as defined in Corollary \ref{diagonals}.  For $n\geq 5$ and $\bd = (d_0,\hdots,d_{n-1})\in\calAC_n$, define:
$$ D_0(\bd) = A(d_0,d_1,\ell_{0,2}(\bd)) + \left[\sum_{i=3}^{n-2} A(\ell_{0,i-1}(\bd),d_i,\ell_{0,i}(\bd))\right] + A(\ell_{0,n-2}(\bd),d_{n-2},d_{n-1}) $$
For each $n\geq 3$, the function $D_0\co\calAC_n\to\mathbb{R}^+$ so-defined is smooth.\end{corollary}

\begin{proof} Since the functions $\ell_{i,j}$ are smooth by Lemma \ref{diagonals}, and the function $A$ from Lemma \ref{triangle area} is smooth on the domain $U$ described there, to establish smoothness of $D_0$ we must only show that each of the different inputs to $A$ that occur in its definition lie in $U$.  This follows from the fact that each is the side length collection of a cyclic triangle: for instance, given $\bd= (d_0,\hdots,d_{n-1})\in\calAC_n$ ($n\geq 5$) and a cyclic $n$-gon $\{x_0,\hdots,x_{n-1}\}$ with side length collection $\bd$, $(\ell_{0,i-1}(\bd),d_i,\ell_{0,i}(\bd))$ is the side length collection of the cyclic triangle $\{x_0,x_{i-1},x_i\}$.

If $(a,b,c)$ is the side length collection of a cyclic triangle then $\sinh(a/2) < \sinh(b/2) + \sinh(c/2)$ by Proposition \ref{existence uniqueness}.  Since the hyperbolic sine has positive first and second derivatives on $(0,\infty)$ it is superadditive there, so $\sinh(b/2)+\sinh(c/2) < \sinh((b+c)/2)$, and we conclude that $a< b+c$.  Reordering and repeating this argument shows $(a,b,c)\in U$.\end{proof}

\section{Geometry (and more calculus)}

As it currently stands following Definition \ref{cyclic order}, for us a ``convex cyclic $n$-gon'' is simply a sequence of points on a hyperbolic circle.  Here we will first show that such a collection is the vertex set of a convex polygon in the classical sense --- a finite intersection of hyperbolic half-spaces --- which is compact, with area measured by the function $D_0$ of Corollary \ref{area function}.

\begin{lemma}\label{polygon}  For $n\geq 3$, a cyclic $n$-gon $\{x_0,\hdots,x_{n-1}\}$ in a hyperbolic circle $C$ is the vertex set of its convex hull, a compact, convex polygon $P$ contained in the disk bounded by $C$.  The edges of $P$ are the geodesic arcs $\gamma_i$ joining $x_{i-1}$ to $x_i$ for each $i>0$, together with $\gamma_0$ joining $x_{n-1}$ to $x_0$.  The area of $P$ is $D_0(d_0,\hdots,d_{n-1})$, where $d_i$ is the length of $\gamma_i$ for each $i$.

Conversely, if the vertex set of a compact, convex polygon $P$ lies in a hyperbolic circle then enumerating it $\{x_0,\hdots,x_{n-1}\}$ so that with the boundary orientation from $P$ an edge points from $x_{i-1}$ to $x_i$ for each $i>0$, and from $x_{n-1}$ to $x_0$, yields a cyclic $n$-gon in the sense of Definition \ref{cyclic order}.\end{lemma}

\begin{proof}  For each $i > 0$ let $\calh_i$ be the half-space bounded by the geodesic through $x_{i-1}$ and $x_i$ such that $\calh_i \cap [x_{i-1},x_i] = \{x_{i-1},x_i\}$.  Define $\calh_0$ analogously so that $\calh_0 \cap[x_{n-1},x_0] = \{x_{n-1},x_0\}$, and let $P = \bigcap_{i=0}^{n-1} \calh_i$.  Since the $\calh_i$ are closed and convex in $\mathbb{H}^2$, so is $P$.

The criterion that $\{x_0,\hdots,x_{n-1}\}$ be cyclically ordered ensures for each $i$ that all $x_j$ other than $x_{i-1}$ and $x_i$ are contained in the interior of $\calh_i$, since $\partial\calh_i$ intersects $C$ only in $\{x_{i-1},x_i\}$.  In particular, $\{x_0,\hdots,x_{n-1}\}\subset P$.  The frontier $\partial P$ of $P$ in $\mathbb{H}^2$ is contained in $\bigcup_{i=0}^{n-1} \partial \calh_i$, and by the above $\gamma_i = \calh_i\cap\partial P$ for each $i$, since $\partial \calh_i$ exits $\calh_{i+1}$ at $x_i$ and $\calh_{i-1}$ at $x_{i-1}$.

We will appeal to \cite{Ratcliffe} for basic results on polygons.  $P$ satisfies the definition of polygon in Section 6.3 there: it is closed, convex and non-empty, and its collection of \textit{sides} (defined in Section 6.2 there) is $\{\gamma_i\}_{i=0}^{n-1}$.  The $\gamma_i$ are thus its edges (``$1$-faces'' in the notation of \cite[\S 6.3]{Ratcliffe}), and the $x_i$ are its vertices, being endpoints of the $\gamma_i$.  $P$ moreover satisfies the compactness criterion of \cite[Theorem 6.3.7]{Ratcliffe}, so by Theorem 6.3.17 there it is the convex hull of its vertex set.  Therefore $\{x_0,\hdots,x_{n-1}\}$ uniquely determines $P$, and since it is contained in the (convex) disk bounded by $C$, so is $P$.

The diagonals from $x_0$ to $x_2,\hdots,x_{n-2}$ divide $P$ into a non-overlapping union of cyclic triangles, whose area sums to the area of $P$.  The diagonal lengths are $\ell_{0,i}(\bd)$ for $i\in\{2,\hdots,n-2\}$, where $\bd = (d_0,\hdots,d_{n-1})$ is the side length collection of $\{x_0,\hdots,x_{n-1}\}$.  Therefore Lemma \ref{triangle area} implies that $D_0(\bd)$ is the area of $P$.

The key observation in showing the converse statement is that given a polygon $P$ inscribed in a circle $C$, any edge $\gamma$ of $P$ bounds a bigon outside $P$ with an arc $\gamma$ of $C$, and the initial and terminal vertices of $\gamma$ in the boundary orientation from $P$ agree with those of $\gamma$ in the counterclockwise orientation on $C$.
\end{proof}

Below we re-characterize the ``centeredness'' condition from Proposition \ref{Cn and BCn} in geometric terms, and we describe an isosceles decomposition which will be useful in analyzing $D_0$.

\begin{proposition}\label{isosceles decomp} For $n\geq 3$, a cyclic $n$-gon $\{x_0,\hdots,x_{n-1}\}$ is centered if and only if the center $v$ of its circumcircle is contained in the interior of its convex hull $P$.  If this is so then $P$ decomposes as the non-overlapping union $\bigcup_{i=0}^{n-1} T_i$, where $T_i$ is the triangle with vertices $v$, $x_i$ and $x_{i-1}$ for $i>0$ and $T_0$ has vertices $v$, $x_0$ and $x_{n-1}$.

If $\{x_0,\hdots,x_{n-1}\}$ is not centered then $P$ has a unique longest side $\gamma_{i_0}$, characterized by the fact that the geodesic containing $\gamma_{i_0}$ has $v$ and $P$ in opposite half-spaces.  In this case $P\cap T_{i_0} = \gamma_{i_0}$, and $P\cup T_{i_0}$ is a convex polygon that decomposes as the non-overlapping union $\bigcup_{i\neq i_0} T_i$.  

If $\{x_0,\hdots,x_{n-1}\}$ is semicyclic then $v$ is the midpoint of $\gamma_{i_0}$, $T_{i_0}=\gamma_{i_0}$, and also $P = \bigcup_{i=0}^{n-1} T_i$.\end{proposition}

\begin{proof}  Let $\{x_0,\hdots,x_{n-1}\}$ have side length collection $(d_0,\hdots,d_{n-1})$.  Lemma \ref{angle function} implies that the triangle $T_i$ defined above has interior angle $\theta(d_i,J)\in (0,\pi]$ at $v$, where $J$ is the radius of the circumcircle $C$.  The edge of $T_i$ opposite $v$, which is the geodesic arc $\gamma_i$ joining $x_i$ to $x_{i-1}$, divides $C$ into two arcs: $[x_{i-1},x_i]$ and $[x_i,x_{i-1}]$ in the notation of Definition \ref{cyclic order}.  If $T_i$ is not degenerate; i.e.~if $\theta(d_i,J)\neq \pi$ then the shorter of these arcs lies on the opposite side of the geodesic $\delta_i$ containing $\gamma_i$ from $v$.

It follows that for each $i$ such that the angle from $x_{i-1}$ to $x_i$ is less than $\pi$, it is $\theta(d_i,J)$, and $v$ lies in the interior of the half-space $\calh_i$ bounded by $\delta_i$ that intersects $[x_{i-1},x_i]$ only in $\{x_{i-1},x_i\}$.  If the angle from $x_{i-1}$ to $x_i$ is greater than $\pi$ then it is $2\pi-\theta(d_i,J)$, and $v$ lies in the interior of the half-space $\calh_i'$ bounded by $\delta_i$ opposite $\calh_i$.  The angle from $x_{i-1}$ to $x_i$ is $\pi$ if and only if $v$ is the midpoint of $\gamma_i$ (the degenerate case mentioned above), so in this case $v$ is in $\calh_i\cap\calh_i'$ but in the interior of neither.

Recall from Proposition \ref{Cn and BCn} that $\{x_0,\hdots,x_{n-1}\}$ is centered if and only if the angle from $x_{i-1}$ to $x_i$ is less than $\pi$ for each $i$.  By the paragraph above this holds if and only if $v$ is in the interior of $P = \bigcap_{i=0}^{n-1} \calh_i$ (compare the proof of Lemma \ref{polygon}). In this case the decomposition of $P$ as a non-overlapping union of the $T_i$ is obtained by simply coning from $v$ to $\partial P$ (which by Lemma \ref{polygon} is the union of the $\gamma_i$).  Again by Proposition \ref{Cn and BCn}, $\{x_0,\hdots,x_{n-1}\}$ is semicyclic if and only if $J=d_{i_0}/2$, where $d_{i_0}$ is maximal among the $d_i$.  In this case $\gamma_{i_0}$ is a diameter of $C$, so $v$ is its midpoint and $T_{i_0} = \gamma_{i_0}$ as claimed, and again coning from $v$ gives $P = \bigcup T_i$.

If the angle from $x_{i_0-1}$ to $x_{i_0}$ is at least $\pi$ for some $i_0$ this $i_0$ is unique, since the angles from $x_{i-1}$ to $x_i$ sum to $2\pi$ ($\{x_0,\hdots,x_{n-1}\}$ being cyclically ordered).  Therefore this angle is $2\pi-\theta(d_{i_0},J)$, for all other $i$ the angle from $x_{i-1}$ to $x_i$ is $\theta(d_i,J)<\pi$, and $\theta(d_{i_0},J) = \sum_{i\neq i_0} \theta(d_i,J)$.  In particular, $\theta(d_{i_0},J) >\theta(d_i,J)$ so $d_{i_0}>d_i$ for all $i\neq i_0$.  In this case $v$ lies in $\calh_{i_0}'$, hence so does $T_{i_0}$, and $T_{i_0}\cap P = \gamma_{i_0}$; but $v$ is in the interior of $\calh_i$ for all $i\neq i_0$. 

Because the angle from $x_{i_0-1}$ to $x_{i_0}$ is at least $\pi$ the diameters of $C$ through $x_{i_0-1}$ and $x_{i_0}$ bound half-spaces $\calh_{i_0-1}'$ and $\calh_{i_0}'$, respectively, which contain all $x_i$ and hence $P$.  It follows as in the proof of Lemma \ref{polygon} that $P\cup T_{i_0} = \calh_{i_0-1}'\cap\calh_{i_0}'\cap\bigcap_{i\neq i_0} \calh_i$ is a convex polygon, and its decomposition as $\bigcup_{i\neq i_0} T_i$ follows by coning from $v$.
\end{proof}

\begin{proposition}\label{smooth D}\SmoothD\end{proposition}

\begin{proof}  We already showed in Corollary \ref{area function} that $D_0$ is smooth, and in Lemma \ref{polygon} that $D_0(d_0,\hdots,d_{n-1})$ is the area of the convex hull of a cyclic $n$-gon with side length collection $(d_0,\hdots,d_{n-1})$ .  Given the decomposition of Proposition \ref{isosceles decomp}, we can rewrite its formula as follows: for a given $\bd = (d_0,\hdots,d_{n-1})\in\calAC_n$,
\begin{align}\label{defect dichot}
 D_0(\bd) = \left\{\begin{array}{ll}
   \sum_{i=0}^{n-1} A(J,J,d_i) & \mbox{if}\ \bd\in\calc_n,\ \mbox{where}\ J = J(\bd) \\
   \left(\sum_{i\neq i_0} A(J,J,d_i)\right) - A(J,J,d_{i_0}) & \mbox{otherwise, where}\ d_{i_0} = \max\{d_i\} \end{array}\right. \end{align}
This follows from Lemma \ref{triangle area} and the fact that each of the triangles $T_i$ of Proposition \ref{isosceles decomp} is isosceles, with two sides of length $J(\bd)$ and one of length $d_i$.  It implies that $D_0$ is symmetric in $(d_0,\hdots,d_{n-1})$.

Comparing the function $\alpha$ of Lemma \ref{triangle area} with $\theta$ from Lemma \ref{angle function}, we note that $\alpha(d,J,J) = \theta(d,J)$ for any $d<2J$; moreover a little hyperbolic trigonometry shows that $\alpha(J,d,J)=\alpha(J,J,d) = \cos^{-1}(\coth J\tanh (d/2))$.  These facts and some trigonometric identities can be used to show the following:
$$ \cos(A(J,J,d_i)/2) = \frac{\cosh^2(d_i/2) +\cosh J}{\cosh(d_i/2)(\cosh J+1)} $$
Taking derivatives and doing some more trigonometry gives:
$$ \frac{\partial A}{\partial d_j}(J,J,d_i) = \left\{\begin{array}{ll}
  \frac{2\sinh J}{\cosh J+1}\frac{\sinh(d_i/2)}{\sqrt{\sinh^2 J - \sinh^2(d_i/2)}} \frac{\partial J}{\partial d_j} & j \neq i \\
  \frac{2\sinh J}{\cosh J+1}\frac{\sinh(d_j/2)}{\sqrt{\sinh^2 J - \sinh^2(d_j/2)}} \frac{\partial J}{\partial d_j} - \frac{\cosh^2(d_j/2) - \cosh J}{\cosh(d_j/2)\sqrt{\sinh^2 J  - \sinh^2(d_j/2)}} & j = i
\end{array}\right. $$
Note that this formula is only defined for $d_i < 2J$.  (This is because the inverse cosine is smooth only on $(-1,1)$.)  We will use it below to compute partial derivatives of $D_0$ at $\bd = (d_0,\hdots,d_{n-1})\in\calAC_n$ but on account of this issue will assume that $J(\bd) < \max\{d_i\}$; i.e.~that $\bd\notin\calBC_n$ (recall Proposition \ref{Cn and BCn}).  Since $\calBC_n$ is a codimension-one submanifold (Proposition \ref{submanifold}), values there are determined by continuity.  

Using symmetricity of $D_0$ we will assume below that $d_0$ is maximal among the $d_i$.  In this case dividing out a factor of $\sqrt{\sinh^2 J - \sinh^2(d_0/2)}$ in the formula of Lemma \ref{J derivs} gives:
$$ \frac{\partial J}{\partial d_j}(\bd) = \frac{1}{2}\frac{\pm\cosh(d_j/2)\tanh J}{\sqrt{\sinh^2 J - \sinh^2(d_j/2)}\left(\frac{\sinh(d_0/2)}{\sqrt{\sinh^2 J - \sinh^2(d_0/2)}}\pm\sum_{i=1}^{n-1} \frac{\sinh(d_i/2)}{\sqrt{\sinh^2 J - \sinh^2(d_i/2)}}\right)} $$
We recall from Lemma \ref{J derivs} that each ``$\pm$'' above should be read as ``$+$'' if $\bd\in\calc_n$ and ``$-$'' if $\bd\in\calAC_n-\calc_n$, except that the numerator is positive in the latter case for $j = 0$. 

We finally come to the derivative computation.  We use the formula above for $D_0$ and consider three cases.  For the first, $\bd\in\calc_n$, we have:\begin{align*}
  \frac{\partial}{\partial d_j} D_0(\bd) & = 
  \frac{2\sinh J}{\cosh J+1}\frac{\partial J}{\partial d_j}\left(\sum_{i=0}^{n-1} \frac{\sinh(d_i/2)}{\sqrt{\sinh^2 J - \sinh^2(d_i/2)}}\right) \\
    & \qquad\qquad\qquad\qquad\qquad\qquad - \frac{\cosh^2(d_j/2) - \cosh J}{\cosh(d_j/2)\sqrt{\sinh^2 J  - \sinh^2(d_j/2)}} \\
    & = \frac{1}{\sqrt{\sinh^2 J - \sinh^2(d_j/2)}}\left[\frac{\sinh^2 J\cosh(d_j/2)}{\cosh J(\cosh J+1)} - \frac{\cosh^2(d_j/2) - \cosh J}{\cosh(d_j/2)}\right] \\
  & = \frac{1}{\sqrt{\sinh^2 J - \sinh^2(d_j/2)}}\left[\frac{\cosh^2 J - \cosh^2(d_j/2)}{\cosh J\cosh(d_j/2)} \right] = \frac{\sqrt{\cosh^2 J - \cosh(d_j/2)}}{\cosh J\cosh(d_j/2)} \end{align*}
For $\bd\in\calAC_n-\calc_n$ we first treat the case $\partial D_0/\partial d_j$ for $j>0$, where:
\begin{align*}
   \frac{\partial}{\partial d_j} D_0(\bd) & =  \frac{2\sinh J}{\cosh J+1}\frac{\partial J}{\partial d_j}\left[\left(\sum_{i=1}^{n-1} \frac{\sinh(d_i/2)}{\sqrt{\sinh^2 J - \sinh^2(d_i/2)}}\right) - \frac{\sinh(d_0/2)}{\sqrt{\sinh^2 J - \sinh^2(d_0/2)}} \right] \\ 
  & \qquad - \frac{\cosh^2(d_j/2) - \cosh J}{\cosh(d_j/2)\sqrt{\sinh^2 J  - \sinh^2(d_j/2)}} 
    = \frac{\sqrt{\cosh^2 J - \cosh(d_j/2)}}{\cosh J\cosh(d_j/2)} \end{align*}
(Intermediate steps parallel the previous computation).  Finally, again for $\bd\in\calAC_n-\calc_n$:\begin{align*}
 \frac{\partial}{\partial d_0} D_0(\bd) &= 
  \frac{2\sinh J}{\cosh J+1}\frac{\partial J}{\partial d_0}\left[\left(\sum_{i=1}^{n-1} \frac{\sinh(d_i/2)}{\sqrt{\sinh^2 J - \sinh^2(d_i/2)}}\right) - \frac{\sinh(d_0/2)}{\sqrt{\sinh^2 J - \sinh^2(d_0/2)}} \right] \\
  &\qquad + \frac{\cosh^2(d_0/2) - \cosh J}{\cosh(d_0/2)\sqrt{\sinh^2 J  - \sinh^2(d_0/2)}} 
    = - \frac{\sqrt{\cosh^2 J - \cosh(d_j/2)}}{\cosh J\cosh(d_0/2)} \end{align*}
This proves the result.\end{proof}

\begin{corollary}\label{monotonicity}\Monotonicity\end{corollary}

\begin{proof}  Given such $(d_0,\hdots,d_{n-1})$ and $(d_0',\hdots,d_{n-1}')$ in $\calc_n\cup\calBC_n$, since $D_0$ is symmetric we may assume $d_i \leq d_i'$ for each $i$; and furthermore we will take $d_0$ maximal among the $d_i$.  We will produce a path $(d_0(t),\hdots,d_{n-1}(t))$ from $(d_0,\hdots,d_{n-1})$ to $(d_0',\hdots,d_{n-1}')$, with its interior in $\calc_n$, that is piecewise-smooth and has each $d_i(t)$ non-decreasing.  The result will thus follow directly from the chain rule and Proposition \ref{smooth D}.

In defining the path we will take $d= \min_i\{d_i\}$ and $D = \max_i \{d_i'\}$.  Then:
$$ d_i(t) = \left\{\begin{array}{lcl} 
  d_i & & d+t \leq d_i \\
  d + t & & d_i \leq d+t \leq d_i' \\
  d_i' & & \mbox{otherwise}  \end{array}\right. \qquad\quad\mbox{for}\ 0\leq t\leq D-d $$
It is clear by inspection that $d_i(t)$ is non-decreasing and piecewise-smooth.  Thus it remains only to check that $(d_0(t),\hdots,d_{n-1}(t))\in\calc_n$.  We break this up into cases.

For $0\leq t \leq d_0 - d$, $d_0 = d_0(t)$ is maximal among the $d_i(t)$ by construction, since $d_0$ is maximal among the $d_i$.  Here we have:
$$ \sum_{i=0}^{n-1} \theta(d_i(t),d_0(t)/2) = \sum_{i=0}^{n-1} \theta(d_i(t),d_0/2) > \sum_{i=0}^{n-1} \theta(d_i,d_0/2) \geq 2\pi $$
The first inequality above follows from the fact that $\theta(d,J)$ (introduced in Lemma \ref{angle function} increases in $d$.  The second follows from the hypothesis that $(d_0,\hdots,d_{n-1})\in\calc_n\cup\calBC_n$ and Proposition \ref{Cn and BCn}, which thus also implies that $(d_0(t),\hdots,d_{n-1}(t))\in\calc_n$.

If $d_0 = D$ then the above case completes the proof, so assume $D > d_0$.  Let $i_0\in\{0,\hdots,n-1\}$ be such that $D = d_{i_0}'$, and let $d_{i_1}'$ be maximal among the $d_i'$ with $i\neq i_0$.

First suppose that $d_0 < d_{i_1}'$.  Then for $d_0-d \leq t \leq d_{i_1}' - d$, at least $d_{i_0}(t)$ and $d_{i_1}(t)$ take the maximum value $d+t$ among the entries $d_i(t)$, so:
$$ \sum_{i = 0}^{n-1} \theta(d_i(t),(d+t)/2) = 2\pi + \sum_{i \neq i_0,i_1} \theta(d_i(t),(d+t)/2) > 2\pi $$
Thus $(d_0(t),\hdots,d_{n-1}(t))\in\calc_n$ for $d_0 - d \leq t\leq d_{i_1}'-d$.

We finally consider the interval $d_{i_1}'-d \leq t < D-d$.  Note that if $d_0 \geq d'_{i_1}$ then this and the interval $0\leq t\leq d_0-d$ cover the domain of the $d_i(t)$, and we can skip the case above.  On the other hand if $d_{i_1}'=d_{i_0}'$ then the previous two cases cover the entire domain and we are done.  Let us therefore assume that $d_{i_1}' < d_{i_0}' = D$.  Then on the interval in question, $d_i(t) = d_i'$ for each $i \neq {i_0}$, and $d_{i_0}(t) = d+t$ is the unique maximal $d_i(t)$.  For $t < D-d$ we have:
$$ \sum_{i =0}^{n-1} \theta(d_i(t),d_{i_0}(t)/2) = \pi + \sum_{i \neq i_0} \theta(d_i',d_{i_0}(t)/2) > \pi+ \sum_{i\neq i_0} \theta(d_i',d_{i_0}'/2) \geq 2\pi $$
Thus $(d_0(t),\hdots,d_{n-1}(t))\in\calc_n$ for $d_{i_1}' - d \leq t< D-d$, and the result is proved.
\end{proof}

\section{To infinity...}\label{to infinity}

\begin{definition}\label{ideal boundary}  Let $\overline{\mathbb{H}}^2$ be the closure of the upper half-plane $\mathbb{H}^2$ in the one-point compactification $\mathbb{C}\cup\{\infty\}$ of $\mathbb{C}$. The \textit{ideal boundary} of $\mathbb{H}^2$ is $\overline{\mathbb{H}}^2-\mathbb{H}^2 = \mathbb{R}\cup\{\infty\}$.\end{definition}

In this section we will show that $\calAC_n$, which parametrizes cyclic hyperbolic $n$-gons, has as its frontier in $(0,\infty)^n$ a space $\calHC_n$ that parametrizes \textit{horocyclic} $n$-gons: those with vertices on a ``horocycle'' (defined below).  We will also describe \textit{horocyclic ideal} $n$-gons, which are natural limits for certain families of cyclic $n$-gons with edge lengths approaching infinity.

\begin{definition}\label{horocycle}  Let $C_{\infty}= \mathbb{R}+i$ and $B_{\infty} = \{z\in\mathbb{C}\,|\, \Im z\geq1\}$, and note that in $\mathbb{C}\cup\{\infty\}$, $\{\infty\} = \overline{C_{\infty}} - C_{\infty}$.  A   \textit{horocycle} of $\mathbb{H}^2$ is a $\mathrm{PSL}(2,\mathbb{R})$-translate of $C_{\infty}$, and its \textit{ideal point} and the \textit{horoball} that it bounds are the corresponding translates of $\infty$ and $B_{\infty}$, respectively.\end{definition}

The action of $\mathrm{PSL}_2(\mathbb{R})$ is transitive on $\mathbb{R}\cup\{\infty\}$, the stabilizer of $\infty$ acts transitively on horizontal lines in $\mathbb{H}^2$ via $\left\{\left(\begin{smallmatrix} r & 0 \\ 0 & 1/r \end{smallmatrix}\right)\right\}$, and the stabilizer $\left\{\left(\begin{smallmatrix} 1 & r \\ 0 & 1 \end{smallmatrix}\right)\right\}$ of $C_{\infty}$ acts transitively on it.

The horocycles are thus the horizontal straight lines and the non-empty intersections with $\mathbb{H}^2$ of circles tangent to $\mathbb{R}$.  The horoball that one bounds is the region above the line in the former case, and inside the circle in the latter.  So it is natural to orient horocycles ``counterclockwise'' by giving them the boundary orientation from their horoballs.

\begin{definition}\label{horocyclic order}  If $C$ is a horocycle with ideal point $v$ then $C\cup\{v\}$ is a circle.  We may thus define the notions of \textit{counterclockwise} and \textit{cyclic order} by analogy with Definition \ref{cyclic order}.

A \textit{horocyclic $n$-gon} is a collection $\{x_0,\hdots,x_{n-1}\}$ of distinct points on a horocycle $C$ with ideal point $v$ that is cyclically ordered on $C\cup\{v\}$.  A \textit{horocyclic ideal $n$-gon} is a cyclically ordered collection $\{x_0,\hdots,x_{n-1}\}$ on some $C\cup\{v\}$ with an \textit{ideal vertex} $x_{i} = v$.

The \textit{side length collection} of a horocyclic $n$-gon $\{x_0,\hdots,x_{n-1}\}$ is $(d_0,\hdots,d_{n-1})$, where $d_i = \mathit{dist}(x_{i-1},x_i)$ for each $i>0$ and $d_0 = \dist(x_0,x_{n-1})$.  It is defined analogously for a horocyclic ideal $n$-gon $\{x_0,\hdots,x_{n-1}\}$, except that $d_i = \infty$ if $x_i$ or $x_{i-1}$ is the ideal vertex.\end{definition}

\begin{proposition}\label{horocyclic parameters}  For $n\geq 3$, $(d_0,\hdots,d_{n-1})\in(0,\infty)^n$ is the side length collection of a horocyclic $n$-gon if and only if $\sinh(d_i/2) = \sum_{j\neq i} \sinh(d_j/2)$ for some $i$; and $(d_0,\hdots,d_{n-1})\in(0,\infty]^n$ is the side length collection of a horocyclic ideal $n$-gon if and only if $d_{i_0} = d_{i_0+1} = \infty$ for a unique $i_0$ (taking $i_0 + 1 = 0$ if $i_0 = n-1$).

Two horocyclic or horocyclic ideal $n$-gons are isometric if and only if their side length collections differ by a cyclic permutation.\end{proposition}

\begin{proof}  Since the isometry group of $\mathbb{H}^2$ acts transitively on horocycles, given a horocyclic $n$-gon $\{x_0,\hdots,x_{n-1}\}$ we may assume it lies on $C_{\infty}$.  The key fact here follows from a short explicit calculation (or an appeal to eg. \cite[Th.~1.2.6(iii)]{Katok}): for $x,y\in C_{\infty}$, $\sinh(\dist(x,y)/2) = |x-y|/2$.

For the unique $i$ such that the counterclockwise arc from $x_{i-1}$ to $x_i$ contains $\{\infty\}$, the real coordinate of $x_i$ is minimal among the $x_j$, and real coordinates increase along the sequence $x_i,\hdots,x_{n-1},x_0,\hdots,x_{i-1}$.  The key fact thus implies that $\sinh(d_i/2) = \sum_{j\neq i}\sinh(d_j/2)$.

On the other hand, given $(d_0,\hdots,d_{n-1})$ such that $\sinh(d_i/2) = \sum_{j\neq i} \sinh(d_j/2)$ one a collection $\{x_0,\hdots,x_{n-1}\}$ on $C_{\infty}$, laid out in the order above, has side length collection $(d_0,\hdots,d_{n-1})$ if the real coordinates of successive points $x_{i-1}$ and $x_i$ differ by $\ell_i = 2\sinh(d_i/2)$.

It is clear that the condition for $(d_0,\hdots,d_{n-1})\in(0,\infty]^n$ to be the side length collection of a horocyclic ideal $n$-gon is necessary.  Existence follows the strategy above: we put $x_{i_0}$ at $\infty$ for the unique $i_0$ such that $d_{i_0} = d_{i_0+1} = \infty$ and arrange $x_{i_0+1},\hdots,x_{n-1},x_0,\hdots,x_{i_0-1}$ on $C_{\infty}$ with real coordinates in increasing order.\end{proof}

\begin{corollary}\label{HCn}  For each $n\geq 3$, set of marked, horocyclic $n$-gons in $\mathbb{H}^2$ is parametrized by the frontier of $\calAC_n$ (defined in Corollary \ref{ACn}) in $(0,\infty)^n$:
$$ \calHC_n = \left\{(d_0,\hdots,d_{n-1})\in(0,\infty)^n\,|\, \sinh(d_{i_0}/2) = \sum_{i\neq i_0} \sinh(d_{i_0}/2),\ \mbox{where}\ d_{i_0}=\max\{d_i\}_{i=0}^{n-1} \right\} $$
It is the orbit of $\mathrm{graph}(h_0) \doteq \{(h_0(\bd),\bd)\,|\,\bd\in(\mathbb{R}^+)^{n-1}\}$ under cyclic permutation of entries, where $h_0(d_1,\hdots,d_{n-1}) = 2\sinh^{-1}\left(\sum_{i=1}^{n-1} \sinh(d_i/2)\right)$ is symmetric, smooth and strictly increasing in each variable.  Moreover, $h_0(\bd) > b_0(\bd)$ for each $\bd\in(\mathbb{R}^+)^{n-1}$, where $b_0$ is as in Proposition \ref{submanifold}, so $\calHC_n$ has a neighborhood disjoint from $\calc_n\cup\calBC_n$ in $(0,\infty)^n$.

The set of marked, horocyclic ideal $n$-gons is parametrized by:
$$ \calHI_n = \left\{(d_0,\hdots,d_{n-1})\in(0,\infty]^n\,|\,d_{i_0} = d_{i_0+1} = \infty\ \mbox{for a unique}\ i_0,\ 0\leq i_0<n \right\} $$
It is the orbit of $\{(\infty,\infty)\}\times\mathbb{R}^{n-2}$ under cyclic permutation of entries.\end{corollary}

\begin{proof}  Most of the above is obvious and/or a consequence of Proposition \ref{horocyclic parameters}, but we will prove that $h_0(\bd) > b_0(\bd)$.  Recall from Proposition \ref{submanifold} that $b_0(\bd)$ is defined by the equation $\sum_{i=1}^{n-1} \theta(d_i,b_0(\bd)/2) = \pi$, for $\theta$ as in Lemma \ref{angle function}.  Plugging into $\ell^n$ from Proposition \ref{omnibus1} gives $\ell^n(b_0(\bd)/2,\bd) = b_0(\bd)$.  The inequality thus follows from the fact that fixing $\bd$, $\ell^n$ increases in $J$ to a limit of $h_0(\bd)$, by Proposition \ref{omnibus2}.\end{proof}

\begin{proposition}\label{radius up_n_down}  For $n\geq 3$, values of the circumcircle radius function $J$ approach infinity on a sequence in $\calAC_n$ approaching $\calHC_n$ or $\calHI_n$.\end{proposition}

\begin{proof} Suppose $\bd = (d_0,\hdots,d_{n-1})\in\calHC_n$ is approached by a sequence in $\calAC_n$.  The unique maximal entry of $\bd$ is $d_{i_0}$ such that $\sinh(d_{i_0}/2) = \sum_{i\neq i_0} \sinh(d_i/2)$, so the $i_0$ entry is also maximal for all but finitely many points in the sequence.  Using symmetricity of $J$ we will assume that $i_0 = 0$, so Lemma \ref{reduce} implies that $d_0 = \ell^n(J)$ for each point in the sequence.

If it were not true that $J(\bd)\to\infty$ then upon passing to a subsequence we could ensure that $J(\bd)\to J_0$ for some real $J_0$.  But $\ell^n$ is continuous as a function of $(d_1,\hdots,d_{n-1},J)$ such that $J \geq 2\max\{d_i\}$ (cf.~ Proposition \ref{omnibus2}), so since the initial coordinate also converges this would imply that $d_0 = \ell^n(J_0,d_1,\hdots,d_{n-1})$ at $\bd\in\calHC_n$.  But this is not possible since $\ell^n$ increases toward its asymptote, which is $d_0 = h_0(d_1,\hdots,d_{n-1})$ (recall Proposition \ref{omnibus2}).

That values of $J$ approach infinity on a sequence approaching $\calHI_n$ is an immediate consequence of the fact that $J(d_0,\hdots,d_{n-1})\geq \max\{d_i\}/2$.\end{proof}

\begin{proposition}\label{horocyclic defect} For $n\geq 3$, the formulas below define a symmetric, continuous extension of $D_0$ to $\calAC_n\cup\calHC_n\cup\calHI_n$.  For $(d_0,\hdots,d_{n-1})\in\calHC_n$ with maximal entry $d_{i_0}$, define:
$$ D_0(d_0,\hdots,d_{n-1}) = (n-2)\pi + 2\left[\sin^{-1}\left(\frac{1}{\cosh(d_{i_0}/2)}\right) - \sum_{i\neq i_0} \sin^{-1}\left(\frac{1}{\cosh(d_i/2)}\right)\right] $$
For $(d_0,\hdots,d_{n-1})\in\calHI_n$ with $d_{i_0} = d_{i_0+1} = \infty$, take:
$$ D_0(d_0,\hdots,d_{n-1}) = (n-2)\pi - 2\sum_{i\neq i_0,i_0+1} \sin^{-1}\left(\frac{1}{\cosh(d_i/2)}\right) $$
Given $(d_0,\hdots,d_{n-1})$ and $(d_0',\hdots,d_{n-1}')$ in $\calHC_n\cup\calHI_n$, if up to a fixed permutation $d_i \leq d_i'$ for each $i$, and $d_i<d_i'$ for some $i$, then $D_0(d_0,\hdots,d_{n-1}) < D_0(d_0',\hdots,d_{n-1}')$.\end{proposition}

\begin{proof}  Both cases of this result follow from:

\begin{claim}\label{area limit}  For sequences $\{a_k\}$, $\{b_k\}$ and $\{c_k\}$ of positive real numbers such that $a_k\to a\in(0,\infty)$ and $b_k \to \infty$ as $k\to\infty$, and $|\sinh(c_k/2)-\sinh(b_k/2)|<\sinh(a_k/2)$ for all $k$, the area function $A$ of Lemma \ref{triangle area} satisfies:
$$ \lim_{k\to\infty} A(a_k,b_k,c_k) = \pi - 2\sin^{-1}\left(\frac{1}{\cosh(a/2)}\right) $$\end{claim}

The claim is a consequence of the definition of $A$ and the following limit computation:\begin{align*}
 & \lim_{k\to\infty} \alpha(a_k,b_k,c_k) = 0 \\
 & \lim_{k\to\infty} \alpha(b_k,c_k,a_k) = \lim_{k\to\infty} \alpha(c_k,a_k,b_k) = \cos^{-1}\left(\frac{\cosh a - 1}{\sinh a} \right) = \sin^{-1}\left(\frac{1}{\cosh(a/2)}\right)\end{align*}
This in turn is a calculus exercise using the definition of $\alpha(a,b,c)$ (again in Lemma \ref{triangle area}).  

In the case that a sequence approaches $(d_0,\hdots,d_{n-1})\in\calHC_n$, by Corollary \ref{HCn} we may assume that all but finitely many terms lie in $\calAC_n-\calc_n$ with maximal $i_0$ entry, where $d_{i_0}$ is maximal among the $d_i$.  We therefore apply the claim to the formula (\ref{defect dichot}), with $b_k = c_k = J$.

For a sequence in $\calAC_n$ approaching $(d_0,\hdots,d_{n-1})\in\calHI_n$, where $d_{i_0} = d_{i_0+1} = \infty$, an arbitrary element is the side length collection of a cyclic polygon $\{x_0,\hdots,x_{n-1}\}$ with longest sides containing $x_{i_0}$.  We divide the convex hull of $\{x_0,\hdots,x_{n-1}\}$ into a non-overlapping union of cyclic triangles using diagonals from $x_{i_0}$, with side lengths given by the diagonal functions $\ell_{i_0,j}$ of Corollary \ref{diagonals}, for $j\neq i_0,i_0+1$.  We compute values of $D_0$ by summing the areas of these triangles.  Since they are cyclic their side lengths satisfy the claim's hypotheses.

The resulting extensions of $D_0$ are clearly symmetric and continuous.  For $(d_0,\hdots,d_{n-1})$ and $(d_0',\hdots,d_{n-1}')\in\calHC_n$ such that $d_i \leq d_i'$ for each $i$, if $d_{i_0}$ is maximal among the $d_i$ and $d_{i_1}'$ is maximal among the $d_i'$ then $d_{i_1}' > d_{i_0}' \geq d_{i_0}$ and $d_{i_0}' \geq d_{i_0} > d_{i_1}$.  Thus exchanging $d_{i_1}'$ with $d_{i_0}'$ does not change the fact that $d_i\leq d_i'$ for all $i$.  Having made the exchange, and assuming without loss of generality that $i_0 = 0$, a little calculus shows that $D_0$ increases along the graph of the straight-line path from $(d_1,\hdots,d_{n-1})$ to $(d_1',\hdots,d_{n-1}')$ in $\mathrm{graph}(h_0)$ (recall Corollary \ref{HCn}).

The corresponding monotonicity property of $D_0$ on $\calHI_n$ follows by direct comparison.  To compare $(d_0,\hdots,d_{n-1})\in\calHC_n$ and $(d_0',\hdots,d_{n-1}')\in\calHI_n$, we note that if $d_i \leq d_i'$ for all $i$ then for $i_0$ such that $d_{i_0}$ is maximal among the $d_i$, $d_{i_0}'=\infty$.\end{proof}

\begin{proposition}\label{horocyclic decomp}  Let $C$ be a horocycle with ideal point $v$, and for some $n\geq 3$ suppose $\{x_0,\hdots,x_{n-1}\}\subset C\cup\{v\}$ is a horocyclic ideal $n$-gon with ideal vertex $x_{i_0} = v$.  For each $i>0$ there is a unique half-space $\calh_i$ bounded by the geodesic through $x_i$ and $x_{i-1}$ such that $x_j\in\overline{\calh}_i$ for all $j$, where $\overline{\calh}_i$ is the closure of $\calh_i$ in $\mathbb{C}\cup\{\infty\}$.  Similarly, there is a unique half-space $\calh_0$ bounded by the geodesic through $x_0$ and $x_{n-1}$ such that $x_j\in\overline{\calh}_0$ for all $j$.

Say $P = \bigcap_{i=0}^{n-1} \calh_i$ is the \mbox{\rm convex hull} of $\{x_0,\hdots,x_{n-1}\}$.  It is a convex polygon contained in the horoball bounded by $C$, with vertex set $\{x_0,\hdots,\widehat{x}_{i_0},\hdots,x_{n-1}\}$ and edge set $\{\gamma_i\}$, where $\gamma_i$ is the geodesic joining $x_i$ to $x_{i-1}$ for each $i$.  $P$ decomposes as the non-overlapping union $\bigcup_{i\neq i_0,i_0+1} T_i$, where $T_i$ is the convex hull of $x_i$, $x_{i-1}$ and $v$ for each such $i$. 

Now suppose $\{x_0,\hdots,x_{n-1}\}\subset C$ is a horocyclic $n$-gon with maximal side length $d_{i_0} = \mathit{dist}(x_{i_0-1},x_{i_0})$.  It is the vertex set of its convex hull, a compact, convex polygon $P$ contained in the horoball bounded by $C$.  The edges of $P$ are the geodesic arcs $\gamma_i$ joining $x_{i-1}$ to $x_i$ for each $i>0$, together with $\gamma_0$ joining $x_{n-1}$ to $x_0$.  Taking $T_i$ as in the previous case for each $i$, $P\cup T_{i_0} = \bigcup_{i\neq i_0} T_i$ is a convex polygon.

In each case above the area of $P$ is given by $D_0(d_0,\hdots,d_{n-1})$ from Proposition \ref{horocyclic defect}, where $d_i = \dist(x_i,x_{i-1})$ for each $i$ (in particular, $d_i = \infty$ if $x_i$ or $x_{i-1}$ is $v$).\end{proposition}

\begin{proof}  The geodesics of the upper half-plane $\mathbb{H}^2$ are the intersections with $\mathbb{H}^2$ of vertical Euclidean straight lines and Euclidean circles centered in $\mathbb{R}$.  It follows from this description that any two distinct points of $\mathbb{H}^2$ are contained in the closure in $\mathbb{C}\cup\{\infty\}$ of a unique geodesic.  (In particular, the closure of a vertical Euclidean straight line contains $\infty$.)

Let us now take $C = C_{\infty}$ and address the case that $\{x_0,\hdots,x_{n-1}\}$ is a horocyclic ideal polyhedron.  If $x_{i_0} = \infty$ then $x_{i_0+1},\hdots,x_{n-1},x_0,\hdots,x_{i_0-1}$ lie on $C_{\infty}$ in order of increasing real part.  The half-space $\calh_{i_0+1}$ described above is therefore the region to the right of the vertical Euclidean straight line through $x_{i_0+1}$, and $\calh_{i_0}$ is to the left of the vertical line through $x_{i_0-1}$.  For each $i\neq i_0$ or $i_0+1$, $\calh_i$ is the region outside the Euclidean circle  through $x_i$ and $x_{i-1}$ ($x_{n-1}$ if $i=0$) that is centered in $\mathbb{R}$.

It is now easy to see that the edges of $P = \bigcap \calh_i$ are the $\gamma_i$ described above, and the vertices consist of all $x_i$ but $x_{i_0}$ (compare the proof of Lemma \ref{polygon}).  $P$ is divided into the $T_i$, $i\neq i_0, i_0+1$ by the vertical lines through the $x_i$, $i\neq i_0, i_0\pm 1$.  This is because $T_i$ is bounded by the vertical lines through $x_i$ and $x_{i-1}$, and the circular arc containing both of these points, for each $i$.

If $\{x_0,\hdots,x_{n-1}\}$ is horocyclic but not horocyclic ideal then the proof that it is the vertex set of its convex hull $P$ follows that of Lemma \ref{polygon}, along with the fact that $P$ is compact and the description of the edge set.  Again taking $C = C_{\infty}$, the key difference between this case and the previous one is that since all $x_i$ lie on $C_{\infty}$, all edges are compact circular arcs.  In particular, if $d_{i_0}$ is maximal among the $d_i$ then $x_{i_0}$ has minimal real coordinate, $x_{i_0-1}$ has maximal real coordinate, and $\calh_{i_0}$ is the region \textit{inside} the Euclidean circle through $x_{i_0}$ and $x_{i_0-1}$ that is centered in $\mathbb{R}$.

Given this fact we note that $P\cup T_{i_0}$ is the convex hull of $\{x_0,\hdots,x_{i_0-1},v,x_{i_0},\hdots,x_{n-1}\}$, a horocyclic ideal $(n+1)$-gon, so its decomposition follows from the previous case.

That the area of either $P$ above is $D_0(d_0,\hdots,d_{n-1})$ follows from their decompositions and the fact that $T_i$ has area $\pi - 2\sin^{-1}(1/\cosh(d_i/2))$ for each $i$.  This in turn follows again from the Gauss--Bonnet formula for the area of hyperbolic triangles, see Theorem 3.5.5 of \cite{Ratcliffe}.  $T_i$ is a ``generalized hyperbolic triangle'' in the terminology of \cite[\S 3.5]{Ratcliffe}, with ideal vertex $v$, so its angle at $v$ is defined to be $0$.  An exercise in Euclidean geometry (recalling from the proof of Proposition \ref{horocyclic parameters} that $\sinh (d_i/2) = |x_i-x_{i-1}|/2$, and noting that the hyperbolic and Euclidean metrics on $\mathbb{H}^2$ are conformal; see Definition \ref{uhp}) establishes that its angles at $x_i$ and $x_{i-1}$ are each $\alpha$ satisfying $\sin \alpha = 1/\cosh(d_i/2)$.  (Compare \cite[2.6.12]{Th_notes}.)
\end{proof}

\section{...and beyond!}\label{beyond}

To this point we have proved that $(d_0,\hdots,d_{n-1})\in(0,\infty)^n$ is the side length collection of a unique cyclic or horocyclic $n$-gon in $\mathbb{H}^2$ if and only if $\sinh(d_i/2) \leq \sum_{j\neq i} \sinh(d_j/2)$ for all $i$.  This condition implies that $d_i < \sum_{j\neq i} d_j$, since the hyperbolic sine has positive first and second derivatives on $(0,\infty)$.  Here we address the remaining case: any $(d_0,\hdots,d_{n-1})\in(0,\infty)^n$ such that $\sinh(d_i/2) > \sum_{j\neq i} \sinh(d_j/2)$ but $d_i\leq \sum_{j\neq i} d_j$ for some $i$ is the side length collection of an ``equidistant'' polygon, with vertices equidistant from a fixed geodesic.  

In the upper half-plane model, such an equidistant locus is the intersection with $\mathbb{H}^2$ of a circle in $\mathbb{C}$ not entirely contained in $\overline{\mathbb{H}}^2$.  It's literally beyond infinity!

\begin{definition}  For $J\geq 0$ the \textit{$J$-equidistant locus} to a geodesic $\gamma$ in $\mathbb{H}^2$ is the collection of points that have distance $J$ from $\gamma$.  For $n\geq 3$, a collection $\{x_0,\hdots,x_{n-1}\}$ of distinct points on a component $C$ of the $J$-equidistant locus to $\gamma$ is an \textit{equidistant} $n$-gon if for some $i_0$ the collection $\{x_{i_0},\hdots,x_{n-1},x_0,\hdots,x_{i_0-1}\}$ is \textit{linearly ordered} on $C$: in the orientation that $C$ inherits as a boundary component of the region it bounds with $\gamma$, the compact arc of $C$ bounded by $x_i$ and $x_{i-1}$ points from $x_{i-1}$ to $x_i$ for each $i\neq i_0$.  (Here we take $0-1$ to be $n-1$.)  We say $x_{i_0}$ is \textit{first} among the $x_i$.

The \textit{collar radius} of $\{x_0,\hdots,x_{n-1}\}$ is the distance from the $x_i$ to $\gamma$.  Its \textit{side length collection} is $(d_0,\hdots,d_{n-1})\in(0,\infty)^n$, where $d_i = \dist(x_{i-1},x_i)$ for $i>0$ and $d_0 = \dist(x_{n-1},x_0)$.\end{definition}

We prove existence and uniqueness of equidistant $n$-gons by a strategy parallel to the one for cyclic $n$-gons.  The role of the circle center is played here by the geodesic $\gamma$, and the angle from $x$ to $y$ by the distance from the orthogonal projection of $x$ to the projection of $y$.

\begin{lemma}\label{Dangle function} If points $x$ and $y$ in $\mathbb{H}^2$ with $\dist(x,y) = d\geq 0$ each have distance $J\geq 0$ from a hyperbolic geodesic $\gamma$, then the orthogonal projections of $x$ and $y$ to $\gamma$ are at distance:
$$ \psi(d,J) = 2\sinh^{-1}(\sinh(d/2)/\cosh J) $$
This is a continuous function on $[0,\infty)^2$, smooth in its interior.  For fixed $d>0$, $\psi(d,J)$ decreases in $J$ on $[0,\infty)$ with $\psi(d,0) = d$, $\lim_{J\to\infty} \psi(d,J) = 0$ and
$$ \frac{\partial}{\partial J}\psi(d,J) = \frac{-2\sinh(d/2)\sinh J}{\cosh J\sqrt{\cosh^2 J+\sinh^2(d/2)}} 
  = -2\tanh J\tanh (\psi(d,J)/2)$$
\end{lemma}

\begin{proof} Let $p$ and $q$ be the orthogonal projections of $x$ and $y$ to $\gamma$, respectively.  The arc joining $x$ to $p$ meets $\gamma$ perpendicularly (hence the term ``orthogonal projection''), as does the arc joining $y$ to $q$, so these arcs and hence also $x$ and $y$ are exchanged by reflection in the perpendicular bisector of the arc of $\gamma$ joining $p$ to $q$.  If $Q$ is the quadrilateral with vertices at $x$, $y$, $p$, and $q$, this perpendicular bisector thus divides $Q$ into isometric quadrilaterals $Q_0$ and $Q_1$ with three right angles each.  

Suppose $Q_0$ contains $x$.  If $\alpha$ is its angle at $x$ then hyperbolic trigonometry gives the following relations between $\alpha$ and the side lengths:\begin{align*}
  & \cos\alpha = \sinh h \sinh(\psi(d,J)/2) &
  & \cosh(d/2) = \frac{\cosh(\psi(d,J)/2)}{\sin\alpha} &
  & \cosh J  = \frac{\cosh h}{\sin\alpha} \end{align*}
Here $h$ is the length of the side of $Q_0$ that lies in the perpendicular bisector of the geodesic from $x$ to $y$.  The left-hand equation above follows from Theorem 3.5.10 of \cite{Ratcliffe} and the other two from Theorem 3.5.7 there.  Some manipulations give the formula for $\psi(d,J)$, and also:
\begin{align}\label{angle measure} \sin\alpha = \frac{\sqrt{\cosh^2 J+\sinh^2(d/2)}}{\cosh J\cosh(d/2)}
\end{align}
The remaining assertions are straightforward.\end{proof}

\begin{lemma}\label{Domnibus1}  If $\{x_0,\hdots,x_{n-1}\}$ is an equidistant $n$-gon with collar radius $J$ and side length collection $(d_0,\hdots,d_{n-1})$, such that $x_{i_0}$ is first among the $x_i$, then:
$$ d_{i_0}= L(J,d_0,\hdots,\widehat{d_{i_0}},\hdots,d_{n-1}) \doteq 2\sinh^{-1}\left[\cosh J \sinh\left(\frac{1}{2}\sum_{i\neq i_0} \psi(d_i,J)\right)\right] $$
Moreover, given any fixed collection of $d_i>0$, for $i\neq i_0$, and $J\geq 0$ there is an equidistant $n$-gon $\{x_0,\hdots,x_{n-1}\}$ satisfying the hypotheses above, and any two such are isometric.\end{lemma}

\begin{proof}  Given $d_i>0$ for $i\neq i_0$, and $J\geq 0$, fix a hyperbolic geodesic $\gamma$, a component $C$ of the $J$-equidistant locus to $\gamma$, and $x_{i_0}\in C$.  Let $p\co\mathbb{H}^2\to \gamma$ be the orthogonal projection, and for each $y\in\gamma$ let $p^{-1}(y)$ be the unique intersection point of $C$ with the preimage of $y$ under $p$.

Let $y_{i_0} = p(x_{i_0})$.  We will recursively produce the $x_i$ for $i\neq i_0$, starting at $i = i_0+1$ and taking $i$ modulo $n$: let $y_i$ be the point on $\gamma$ at distance $\psi(d_i,J)$ from $y_{i-1}$ (mod $n$) with the property that the compact arc bounded by $y_i$ and $y_{i-1}$ points toward $y_{i-1}$ in the boundary orientation that $\gamma$ inherits from the region it bounds with $C$.  Note that the restriction of $p$ to $C\to\gamma$ reverses this orientation, so taking $x_i = p^{-1}(y_i)$ for each $i$ we obtain a cyclic $n$-gon on $C$ with $x_{i_0}$ first among the $x_i$.  

Lemma \ref{Dangle function} implies that $\dist(x_i,x_{i-1})=d_i$ for each $i\neq i_0$.  By construction the distance from $y_{i_0}$ to $y_{i_0-1}$ is $\sum_{i\neq i_0} \psi(d_i,J)$, so another application of Lemma \ref{Dangle function} implies that $\dist(x_{i_0},x_{i_0-1})$ satisfies the formula for $d_{i_0}$ described above.  It follows that $\{x_0,\hdots,x_{n-1}\}$ is an equidistant $n$-gon with side length collection $(d_0,\hdots,d_{n-1})$.  By construction it has collar radius $J$ and $x_{i_0}$ first among the $x_i$.

Another equidistant $n$-gon with the properties of $\{x_0,\hdots,x_{n-1}\}$, on a component of the equidistant locus to a hyperbolic geodesic $\gamma'$, can be taken to $\{x_0,\hdots,x_{n-1}\}$ by the following sequence of isometries: take $\gamma'$ to $\gamma$, apply an order-two rotation around a point of $\gamma$ if necessary (to exchange components of the $J$-equidistant locus), then translate in $\gamma$.\end{proof}

\begin{lemma}\label{Domnibus2}  For any $n\geq 3$, the function $L$ defined in Lemma \ref{Domnibus1} is continuous on $[0,\infty)^n$ and smooth on $(0,\infty)^n$.  Fixing $(d_1,\hdots,d_{n-1})$ and taking the restriction of $L$ to $[0,\infty)\times\{(d_1,\hdots,d_{n-1})\}$ as a function of $J$, we have $\frac{\partial}{\partial J} L(J) < 0$ for all $J>0$, $L(0) = \sum_{i=1}^{n-1} d_i$ and $\lim_{J\to\infty} L(J) = 2\sinh^{-1}\left(\sum_{i=1}^{n-1} \sinh(d_i/2)\right)$.\end{lemma}

\begin{proof}  The value of $L(0)$ follows directly from its definition, as does the computation below:\begin{align}\label{collar deriv}
  \frac{1}{2}\cosh (L(J)/2) \frac{\partial L}{\partial J} & = \sinh J \left[\sinh\left(\frac{1}{2}\sum_{i=1}^{n-1} \psi(d_i,J)\right) \right. \\
  &\nonumber\qquad - \left.\cosh\left(\frac{1}{2}\sum_{i=1}^{n-1} \psi(d_i,J)\right)\sum_{i=1}^{n-1} \tanh(\psi(d_i,J)/2) \right]\end{align}
(We also appealed to Lemma \ref{Dangle function} for the derivative of $\psi$ with respect to $J$.)  That $\partial L/\partial J < 0$ now follows from the fact that $\sum\tanh x_i > \tanh\left(\sum x_i\right)$.  This follows by induction along the lines of Claim \ref{add tans} in the proof of Proposition \ref{omnibus2}, using the ``angle addition'' identity for hyperbolic tangent: $\tanh(x+y) = (\tanh x+\tanh y)/(1+\tanh x\tanh y)$.

For the limit as $J\to\infty$ we note that the quantity in brackets in the definition of $L$ is of the form $\infty\cdot 0$, so rewriting $\cosh J$ as $\frac{1}{1/\cosh J}$, applying l'H\^opital's rule and simplifying gives:\begin{align*}
  \lim_{J\to\infty} \cosh J\sinh\left(\frac{1}{2}\sum_{i=1}^{n-1} \psi(d_i,J)\right) & = \lim_{J\to\infty} \cosh\left(\frac{1}{2}\sum_{i=1}^{n-1} \psi(d_i,J)\right)\sum_{i=1}^{n-1} \frac{\sinh(d_i/2)\cosh J}{\sqrt{\cosh^2 J+\sinh^2(d_i/2)}} \\
  & = \sum_{i=1}^{n-1} \sinh(d_i/2) \end{align*}
It follows that the limit of $L$ is as described.\end{proof}

\begin{proposition}\label{equidistant existence uniqueness}  The side length collection $(d_0,\hdots,d_{n-1})$ of an equidistant $n$-gon satisfies $d_{i_0} \leq \sum_{i\neq i_0} d_i$ but $\sinh(d_{i_0}/2)>\sum_{i\neq i_0} \sinh(d_i/2)$, where $d_{i_0}$ is maximal among the $d_i$.  For $n\geq 3$, any $(d_0,\hdots,d_{n-1})\in(0,\infty)^n$ satisfying the above inequalities is the side length collection of a unique equidistant $n$-gon up to isometry.  In particular the collar radius is uniquely determined by, and moreover a symmetric function of, $(d_0,\hdots,d_{n-1})$.  Two equidistant $n$-gons are isometric if and only if their side length collections differ by a cyclic permutation.\end{proposition}

\begin{proof}  For an equidistant polygon $\{x_0,\hdots,x_{n-1}\}$ with side length collection $(d_0,\hdots,d_{n-1})$ and collar radius $J$, if $x_{i_0}$ is first among the $x_i$ then Lemma \ref{Domnibus1} implies that $d_{i_0} = L_{i_0}(J) \doteq 2\sinh^{-1}\left[\cosh J \sinh\left(\frac{1}{2}\sum_{i\neq i_0} \psi(d_i,J)\right)\right]$.  By Lemma \ref{Domnibus2}, $L_{i_0}$ is continuous and strictly decreasing on $[0,\infty)$, with $L_{i_0}(J) \leq \sum_{i\neq i_0} d_i$ and $\sinh(L_{i_0}(J)/2) > \sum_{i\neq i_0} \sinh(d_i/2)$ for all $J$.  The same inequalities therefore hold for $d_{i_0}$, which in particular is maximal among the $d_i$. 

On the other hand, given $(d_0,\hdots,d_{n-1})$ such that the maximal entry $d_{i_0}$ satisfies $d_{i_0} \leq \sum_{i\neq i_0} d_i$ but $\sinh(d_{i_0}/2)>\sum_{i\neq i_0} \sinh(d_i/2)$, since $L_{i_0}$ as defined above is continuous on $[0,\infty)$ there exists some $J\geq 0$ such that $d_{i_0} = L_{i_0}(J)$.  Therefore by Lemma \ref{Domnibus1} there is an equidistant $n$-gon $B$ with side length collection $(d_0,\hdots,d_{n-1})$ and collar radius $J$.  Moreover, since $L_{i_0}$ is strictly decreasing this $J$ is uniquely determined by $(d_0,\hdots,d_{n-1})$, so by Lemma \ref{Domnibus1} again, $P$ is unique up to isometry.

We note that cyclically relabeling the vertices of an equidistant polygon produces an isometric (by the identity map) equidistant polygon whose side length collection is obtained from the original by the same cyclic relabeling.  On the other hand it is clear that an isometry of equidistant polygons takes the side length collection of one to a cyclic permutation of the side length collection of the other.

We note that $L(J)$ is symmetric in $(d_1,\hdots,d_{n-1})$.  This is clear by inspecting its definition in Lemma \ref{Domnibus2}.  Since the collar radius $J(d_0,\hdots,d_{n-1})$ of an equidistant $n$-gon with side length collection $(d_0,\hdots,d_{n-1})$ is determined by the equation $d_0 = L(J,d_0,\hdots,\widehat{d}_{i_0},\hdots,d_{n-1})$ it is therefore invariant under any permutation of the entries fixing the $i_0$ place.  But we already showed that cyclically permuting side lengths yields isometric equidistant polygons; thus with identical collar radii.  Symmetricity of collar radius follows.\end{proof}

The following corollary is immediate.  Here we use ``marked'' as in Corollary \ref{ACn}.

\begin{corollary}\label{equidistant parameters}  For $n\geq3$, marked equidistant $n$-gons in $\mathbb{H}^2$ are parametrized up to isometry by:
$$ \calE_n = \left\{ (d_0,\hdots,d_n)\in(0,\infty)^n\,|\, \sinh(d_{i_0}/2) > \sum_{i\neq i_0} \sinh(d_i/2)\ \mbox{but}\ d_{i_0} \leq \sum_{i\neq i_0} d_i,\ \mbox{for some}\ i_0 \right\} $$
The topological frontier of $\calE_n$ in $(0,\infty)^n$ is $\calHC_n\sqcup\{(d_0,\hdots,d_{n-1})\,|\,d_i = \sum_{j\neq i} d_j\ \mbox{for some}\ i\}$.
\end{corollary}

\begin{proposition}\label{equi smooth J}  The function $J\co\calE_n\to\mathbb{R}^n$ that records collar radius of equidistant $n$-gons is symmetric and continuous on $\calE_n$ and smooth on its interior.  If $d_{i_0}$ is maximal among the $d_i$ then:\begin{align*}
  \frac{\partial J}{\partial d_i}(\bd) = \frac{1}{2}\frac{\pm \cosh(d_i/2)\coth J \sqrt{\frac{\cosh^2 J + \sinh^2(d_{i_0}/2)}{\cosh^2 J + \sinh^2(d_i/2)}}}{\sinh(d_{i_0}/2) - \sum_{j\neq i_0}\sinh(d_j/2)\sqrt{\frac{\cosh^2 J + \sinh^2(d_{i_0}/2)}{\cosh^2 J + \sinh^2(d_j/2)}}}  \end{align*}
Here the ``$\pm$'' is ``$+$'' for $i = i_0$ and ``$-$'' otherwise.  Values of $J$ approach infinity on sequences approaching $\calHC_n$, and $J(d_0,\hdots,d_{n-1})=0$ if and only if $d_{i_0} = \sum_{i\neq i_0} d_i$ for some $i$.\end{proposition}

\begin{proof}  By Lemma \ref{Domnibus1}, $J$ satisfies the equation $d_{i_0} = L(J,d_0,\hdots,\widehat{d}_{i_0},\hdots,d_{n-1})$.  Applying Lemma \ref{Domnibus2} and the implicit function theorem we find that $J$ is smooth on $\calE_n$, with:\begin{align*}
  & \frac{\partial J}{\partial d_{i_0}} = \frac{1}{\partial L/\partial J} &
  & \frac{\partial J}{\partial d_i} = -\frac{\partial L/\partial d_i}{\partial L/\partial J}\quad (i\neq i_0) \end{align*}
(Compare the proof of Lemma \ref{J derivs}.)  Keeping in mind that $d_{i_0} = L$ at the point in question, the computation from (\ref{collar deriv}) yields:
$$ \frac{\partial L}{\partial J} = 2\begin{array}{l} \frac{\sinh J\sqrt{\cosh^2 J + \sinh^2(d_{i_0}/2)}}{\cosh J\cosh(d_{i_0}/2)}\left[\frac{\sinh(d_{i_0}/2)}{\sqrt{\cosh^2 J + \sinh^2(d_{i_0}/2)}} - \sum_{i\neq i_0} \frac{\sinh(d_i/2)}{\sqrt{\cosh^2 J+\sinh^2(d_i/2)}} \right] \end{array}$$
Another computation gives:
$$ \frac{\partial L}{\partial d_i} = \frac{\cosh(d_i/2)}{\cosh(d_{i_0}/2)}\sqrt{\frac{\cosh^2 J+\sinh^2(d_{i_0}/2)}{\cosh^2 J+\sinh^2(d_i/2)}} $$
The derivative computation described above follows.  That $J(d_0,\hdots,d_{n-1}) = 0$ if and only if $d_{i_0} = \sum_{i\neq i_0} d_{i_0}$ for some $i_0$ follows from Lemma \ref{Domnibus2}, since $L$ satisfies $L(0) = \sum_{i\neq i_0} d_i$ and $\frac{\partial}{\partial J} L < 0$ on $(0,\infty)$, with $d_{i_0} = L(J,d_0,\hdots,\widehat{d}_{i_0},\hdots,d_{n-1})$ by Lemma \ref{Domnibus1}.

The argument that $J\to\infty$ approaching $\calHC_n$ is essentially identical to the proof of Lemma \ref{radius up_n_down}.  It again uses Lemma \ref{Domnibus2}, since $\lim_{J\to\infty} L(J) = 2\sinh^{-1}\left(\sum_{i\neq i_0} \sinh(d_i/2)\right)$.
\end{proof}

\begin{proposition}\label{equi poly}  For $n\geq 3$, an equidistant $n$-gon $\{x_0,\hdots,x_{n-1}\}$ in a component $C$ of the $J$-equidistant locus to a geodesic $\gamma$ in $\mathbb{H}^2$ is the vertex set of its convex hull, a compact, convex polygon $P$ contained region between $C$ and $\gamma$.  The edges of $P$ are the geodesic arcs $\gamma_i$ joining $x_{i-1}$ to $x_i$ for each $i>0$, together with $\gamma_0$ joining $x_{n-1}$ to $x_0$.

For each $i$ let $Q_i$ be the quadrilateral with vertices at $x_{i-1}$, $x_i$ (if $i=0$, at $x_0$ and $x_{n-1}$) and their projections to $\gamma$.  Then if $x_{i_0}$ is first among the $x_i$, $P\cup Q_{i_0}$ decomposes as the non-overlapping union $\bigcup_{i\neq i_0} Q_i$.\end{proposition}

\begin{proof} The description of $P$ follows as in Lemma \ref{polygon}.  The key fact is that for any $x$ and $y$ on $C$, the geodesic through $x$ and $y$ intersects $C$ in $\{x,y\}$; and it intersects the region bounded by $C$ and $\gamma$ in the geodesic arc joining $x$ to $y$.  In particular, this region is convex. This in turn follows from convexity of the hyperbolic metric (see eg.~the discussion at the beginning of \cite[Chapter II.2]{BrH}; in particular Proposition 2.2 there).

The decomposition of $P\cup Q_{i_0}$ follows similarly to the non-centered case of Proposition \ref{isosceles decomp}.\end{proof}

\begin{proposition}\label{equi defe deri} The function $D_0\co\calE_n\to\mathbb{R}^+$ that records area of equidistant $n$-gons is symmetric and continuous on $\calE_n$ and smooth on its interior.  It satisfies:
$$ \frac{\partial D_0}{\partial d_i} = \left\{\begin{array}{rl}
  -\sqrt{\frac{1}{\sinh^2 J}+\frac{1}{\cosh^2(d_{i_0}/2)}} & i = i_0 \\
  \sqrt{\frac{1}{\sinh^2 J}+\frac{1}{\cosh^2(d_{i}/2)}} & \mbox{otherwise}, \end{array}\right. $$
where $d_{i_0}$ is maximal among the $d_i$.  $D_0(d_0,\hdots,d_{n-1}) = 0$ if and only if $d_{i_0} = \sum_{i\neq i_0} d_i$ for some $i_0$, and $D_0$ extends continuously to $\calHC_n$ by the formula of Proposition \ref{horocyclic defect}.\end{proposition}

\begin{proof}  The proof follows our strategy from the cyclic case.  Let $\{x_0,\hdots,x_{n-1}\}$ be an equidistant polygon with side length collection $(d_0,\hdots,d_{n-1})$, and suppose $x_{i_0}$ is first among the $x_i$.  For each $i$ let $Q_i$ be the quadrilateral defined in Proposition \ref{equi poly}.  As we observed in the proof of Lemma \ref{Dangle function}, $Q_i$ admits a reflection exchanging $x_{i-1}$ with $x_i$, so it has identical angles there.  Call this angle $\alpha_i$.  It is determined by the formula (\ref{angle measure}) with $d_i$ substituted for $d$.  Note that $\alpha_i$ is continuous in $d$ and $J$ on $[0,\infty)^2$ and smooth on $(0,\infty)^2$.  Since $Q_i$ has right angles at its vertices on $\gamma$, it has area $\pi - 2\alpha_i$.

Assuming that $J$ depends on $d$, a computation gives:
$$ \frac{\partial}{\partial d_i}\mathrm{area}(Q_i) = \frac{\frac{\sinh J}{\cosh(d_i/2)}+2\frac{\sinh(d_i/2)}{\cosh J}\frac{\partial J}{\partial d_i}}{\sqrt{\cosh^2 J + \sinh^2(d_i/2)}} $$
For $j\neq i$, the derivative of $Q_i$ with respect to $d_j$ is identical to the above except that the numerator lacks the term $\sinh J/\cosh(d_i/2)$, and there $\partial J/\partial d_j$ replaces $\partial J/\partial d_i$.  Here $J = J(d_0,\hdots,d_{n-1})$ is the collar radius of $\{x_0,\hdots,x_{n-1}\}$.  

The decomposition of Proposition \ref{equi poly} gives $D_0(d_0,\hdots,d_{n-1}) = \left(\sum_{i\neq i_0} \mathrm{area}(Q_i)\right) - \mathrm{area}(Q_{i_0})$.  By Proposition \ref{equi smooth J},  $J(d_0,\hdots,d_{n-1}) = 0$ if and only if $d_{i_0} = \sum_{i\neq i_0} d_i$ for some $i_0$, so it follows from that result and properties of $\mathrm{area}(Q_i)$ described above that $D_0$ is smooth on the interior of $\calE_n$.  Taking derivatives now yields:
$$ \frac{\partial D_0}{\partial d_j} = \left\{\begin{array}{rl}
  \frac{2}{\cosh J}\frac{\partial J}{\partial d_{i_0}}\left[\left(\sum_{i\neq i_0} \frac{\sinh(d_i/2)}{\sqrt{\cosh^2 J + \sinh^2(d_i/2)}}\right) - \frac{\sinh(d_{i_0}/2)}{\sqrt{\cosh^2 J+\sinh^2(d_{i_0}/2)}}\right] \qquad\quad& \\  - \frac{\sinh J}{\cosh(d_{i_0}/2)\sqrt{\cosh^2 J+\sinh^2(d_{i_0}/2)}} & j = i_0 \\
  \frac{2}{\cosh J}\frac{\partial J}{\partial d_{j}}\left[\left(\sum_{i\neq i_0} \frac{\sinh(d_i/2)}{\sqrt{\cosh^2 J + \sinh^2(d_i/2)}}\right) - \frac{\sinh(d_{i_0}/2)}{\sqrt{\cosh^2 J+\sinh^2(d_{i_0}/2)}}\right] \qquad\quad & \\  + \frac{\sinh J}{\cosh(d_{j}/2)\sqrt{\cosh^2 J+\sinh^2(d_j/2)}} & j \neq i_0 \end{array}\right. $$
Substituting for $\partial J/\partial d_j$ and simplifying yields the formulas claimed, noting in particular that $\cosh^2(d_j/2)+\sinh^2 J = \sinh^2(d_j/2)+\cosh^2 J$ by trading a sum with one from one term to the other.

If $(d_0,\hdots,d_{n-1})\in\calE_n$ has $d_{i_0} = \sum_{i\neq i_0} d_i$ for some $i$ then $J(d_0,\hdots,d_{n-1}) = 0$ by Proposition \ref{equi smooth J}, so for each $i$, $\alpha_i = \pi/2$ by the formula (\ref{angle measure}).  It follows that $Q_i$ as described in Proposition \ref{equi poly} has area $0$, so by that result $D_0(d_0,\hdots,d_{n-1}) = 0$ as well.  By the results above $D_0$ is strictly decreasing along the ray 
$$(d_0,\hdots,d_{i_0-1})\times(h_0,\sum_{i\neq i_0} d_i]\times(d_{i_0+1},\hdots,d_{n-1}),$$
so it takes the value $0$ only where $d_{i_0} = \sum_{i\neq i_0} d_i$.

Finally, since values of $J$ limit to infinity on sequences approaching $\calHC_n$ by Proposition \ref{equi smooth J}, values of the $\alpha_i$ determined by (\ref{angle measure}) limit to $\sin^{-1}(\tanh(d_i/2)) = \cos^{-1}(1/\cosh(d_i/2))$.\end{proof}

\begin{corollary}\label{vertical ray} For $\bd = (d_1,\hdots,d_{n-1})\in(0,\infty)^{n-1}$ ($n\geq 3$) and $d_0\geq D = \max\{d_i\}$:
$$(d_0,\bd)\in\left\{\begin{array}{ll} 
  \calc_n & \mbox{for}\ D\leq d_0< b_0(\bd) \\
  \calBC_n & \mbox{for}\ d_0 = b_0(\bd)\ \mbox{(defined in Proposition \ref{submanifold})} \\ 
  \calAC_n - (\calc_n\cup\calBC_n) & \mbox{for}\ b_0(\bd) < d_0 < h_0(\bd) \\
  \calHC_n & \mbox{for}\ d_0 = h_0(\bd)\ \mbox{(defined in Corollary \ref{HCn})} \\
  \calE_n & \mbox{for}\ h_0(\bd) < d_0 \leq \sum_{i=1}^{n-1} d_i \end{array}\right. $$
Moreover, as a function of $d_0$, $D_0(d_0,\bd)$ is continuous, strictly increasing on $[D,b_0(\bd)]$ and strictly decreasing on $\left[b_0(\bd),\sum_{i=1}^{n-1} d_i\right]$.\end{corollary}

\begin{proof}  The breakdown of where points on the ray $[D,\sum d_i]$ lie follows from Proposition \ref{Cn and BCn} (for $\calc_n$ and $\calBC_n$), Corollaries \ref{ACn} (for $\calAC_n$), \ref{HCn} (for $\calHC_n$), and \ref{equidistant parameters} (for $\calE_n$).  Continuity follows from Propositions \ref{horocyclic defect} and \ref{equi defe deri}, monotonicity from Propositions \ref{smooth D} and \ref{equi defe deri}.\end{proof}

\begin{corollary}\label{max area}\MaxArea\end{corollary}

\begin{proof}  The minimal possible final side length of a hyperbolic $n$-gon is $D_0 = \max\{0,d_{i_0}-\sum_{i\neq i_0} d_i\}$, where $d_{i_0}$ is maximal among the $d_i$.  Let us take $\bd = (d_1,\hdots,d_{n-1})\in(0,\infty)^{n-1})$, and consider the values of the area function $D_0$ at points $(d_0,\bd)$ for $d_0 \in \left[D_0,\sum_{i=1}^{n-1} d_i\right]$.  These are maximized at $d_0 = b_0(\bd)$, since by Corollary \ref{vertical ray} this point maximizes the value of $D_0(d_0,\bd)$ for $d_0 \geq D = \max\{d_i\}$, and Propositions \ref{smooth D} and \ref{equi defe deri} imply that $D_0(d_0,\bd)$ increases for $D_0 \leq d_0 \leq D$.

The point $(b_0(\bd),\bd)\in\calBC_n$ represents a semicyclic $n$-gon with $n-1$ sides of length $d_1,\hdots,d_{n-1}$ (see Corollary \ref{vertical ray} above).  Its area maximizes areas of all cyclic, horocyclic or equidistant $n$-gons with $n-1$ sides of length $d_1,\hdots,d_{n-1}$, since $D_0$ is symmetric.  But by Schlenker's Theorem C \cite[pp.~2159--2160]{Schlenker}, the area of all hyperbolic $n$-gons with a fixed side length collection is maximized by the cyclic, horocyclic, or equidistant $n$-gon with that side length collection.  Our result follows.\end{proof}

\section{Degenerations}\label{degenerations}

This section describes the closure $\overline{\calAC}_n$ of $\calAC_n$ in $[0,\infty)^n$.  By Corollary \ref{HCn}, $\calAC_n\cup\calHC_n$ is the closure of $\calAC_n$ in $(0,\infty)^n$.  It thus remains to describe limits of Cauchy sequences in $\calAC_n$ that have some entries approaching $0$.  The limit of most such sequences lie in copies of $\calAC_m$, for $m<n$:

\begin{lemma}\label{shrinking limit}  Fix $n\geq 3$.  For each $m$ with $0<m<n$, let $\cali_{m,n}$ be the collection of $I = (i_0,\hdots,i_{m-1})\in\mathbb{N}^m$ such that $0\leq i_1<i_2<\hdots<i_m\leq n-1$.  For such an $m$ and $I$, define $\phi_I\co\mathbb{R}^m\to\mathbb{R}^n$ by $\phi_I(x_0,\hdots,x_{m-1}) = (y_0,\hdots,y_{n-1})$, where for $0\leq i\leq n-1$,
$$ y_i = \left\{\begin{array}{ll} x_{i_j} & \mbox{if}\ i = i_j\ \mbox{for some}\  0\leq j\leq m-1\\ 0 & \mbox{otherwise}\end{array}\right.$$
Taking $\Delta = \{(r,r)\,|\,r\geq 0\}$, the closure $\overline{\calAC}_n$ of $\calAC_n$ in $[0,\infty)^n$ is:
$$\overline{\calAC}_n = \calAC_n\cup\calHC_n\cup\left(\bigcup_{3\leq m<n,\ I\in\cali_{m,n}} \phi_I(\calAC_m\cup\calHC_m)\right)\cup\left(\bigcup_{I\in\cali_{2,n}} \phi_I(\Delta)\right)$$
Moreover, each such point other than $(0,\hdots,0)$ lies in a unique set above.  The closure of $\calc_n$ has an entirely analogous description, but with each instance of ``$\calAC$'' above replaced by ``$\calc$'' and each instance of ``$\calHC$'' by ``$\calBC$''.\end{lemma}

\begin{proof}  For a sequence in $\calAC_n$ converging in $[0,\infty)^n$ to a point $(d_0,\hdots,d_{n-1})\notin\calAC_n\cup\calHC_n$ at least one entry $d_i$ is equal to $0$.  Note that if all $d_i=0$ then trivially $(d_0,\hdots,d_{n-1})\in\phi_I(\Delta)$ for any $I\in\cali_{2,n}$, so we may assume this does not hold.  Let $d_{i_0},d_{i_1},\hdots,d_{i_m}$ be the non-zero entries, ordered so that $i_0<i_1<\hdots<i_m$, and let $I = (i_0,\hdots,i_m)\in\cali_{m,n}$.

By construction, $(d_0,\hdots,d_{n-1})\in\phi_I((\mathbb{R}^+)^m)$, and it is not hard to see that $m$ and $I$ are unique with this property.  The key observation is that for each $j$ between $0$ and $m-1$, by definition of $\calAC_n$ (recall Corollary \ref{ACn}) and continuity of hyperbolic sine:
$$ \sinh(d_{i_j}/2)\leq \sum_{k\neq j} \sinh(d_{i_k}/2) $$
It follows in short order that $m>1$, if $m=2$ then $d_{i_0} = d_{i_1}$ (i.e.~$(d_0,\hdots,d_{n-1})\in\phi_I(\Delta)$), and if $m\geq 3$ then $(d_0,\hdots,d_{n-1})\in\phi_I(\calAC_m\cup\calHC_m)$.  Moreover, $(d_0,\hdots,d_{n-1})$ lies in $\phi_I(\calHC_m)$ if and only if equality holds for some $j$ above (recall Corollary \ref{HCn}).

A similar argument using the definition of $\calc_n$ (see Proposition \ref{Cn and BCn}) shows that a limit point of $\calc_n$ outside $\calc_n$ is in $\phi_I(\Delta)$ or $\phi_I(\calc_m\cup\calBC_m)$ for some $I$ in $\cali_{2,n}$ or $\cali_{m,n}$ for $m\geq 3$, respectively.

It is finally not difficult to show, for each $m\geq 3$ and $I\in\cali_{m,n}$, that each point of $\phi_I(\calAC_m)$ (or $\phi_I(\calc_m)$) is approached by a sequence in $\calAC_n$ (respectively, $\calc_n$); and moreover that each point of $\phi_I(\Delta)$ is approached by a sequence in $\calc_n$ for any $I\in\cali_{2,n}$.\end{proof}

The parametrizations of $\calBC_n$ and $\calHC_n$ extend similarly to the closure of $(\mathbb{R}^+)^{n-1}$ in $\mathbb{R}^{n-1}$.

\begin{lemma}\label{shrinking parameters}  For $n\geq 3$, the functions $b_0$ and $h_0$ on $(\mathbb{R}^+)^{n-1}$ from Propositions \ref{submanifold} and \ref{HCn}, respectively, extend continuously to $[0,\infty)^{n-1}$.  For $3\leq m < n$ and $I \in \cali_{m-1,n-1}$ (as in Lemma \ref{shrinking limit}), $b_0 = b_0\circ \phi_I$ and $h_0 = h_0\circ\phi_I$ on $(\mathbb{R}^+)^{m-1}$.\end{lemma}

\begin{proof}  For $3\leq m <n$ and $\bd=(d_1,\hdots,d_{m-1})\in(\mathbb{R}^+)^{m-1}$, the equation $\sum_{i=1}^{m-1} \theta(d_i,J) = \pi$ is uniquely solved by $J=\frac{1}{2} b_0(\bd) >\max\{d_i\}/2$.  Since $\theta(d,J)$ decreases in $J$ for any $d>0$ it follows that $\sum_{i=1}^{n-1} \theta(d_i,J_0) > \pi > \sum_{i=1}^{m-1} \theta(d_i,J_1)$ for any  $J_0$ and $J_1$ with $\max\{d_i\}/2 < J_0 < \frac{1}{2}b_0(\bd) < J_1$.  Fixing $I\in \cali_{m-1,n-1}$, it follows that values of the corresponding sum on a fixed sequence in $(\mathbb{R}^+)^{n-1}$ that approaches $\phi_I(\bd)$ greater than $\pi$ at all $J\leq J_0$ and less than $\pi$ for $J\geq J_1$.  Therefore values of $b_0$ on this sequence are eventually contained in $(2J_0,2J_1)$.

The above implies that $b_0$ on $(\mathbb{R}^+)^{n-1}$ extends continuously by $b_0\circ\phi_I$ at $\phi_I(\bd)$.  The proof for $h_0$ is similar.  Let us note also that $b_0(\bd) \to 0$ and $h_0(\bd)\to 0$ if $\bd\to 0$.  This is because the value of either function is bounded above by $(n-1)\max\{d_i\}$ at $(d_1,\hdots,d_{n-1})$, being a side length of an $n$-gon whose other side lengths are the $d_i$.

It remains to consider limits with exactly one nonzero entry.  Applying the definitions directly in this case implies that the extensions must satisfy:
$$ b_0(0,\hdots,0,d,0,\hdots,0) = d = h_0(0,\hdots,0,d,0,\hdots,0) $$
Arguing as above establishes that the so-defined extensions are continuous.
\end{proof}

The radius function $J$ does not extend continuously to $\overline{\calAC}_n$ since it blows up near $\calHC_n$.  Moreover, it is not hard to show that, say $\lim_{x\to 0} h_0 (d,x) = d$, and in consequence that $(d,d,0)$ (or indeed any point of any $\phi_I(\Delta)$) is a limit of horocyclic polygons.  Because $\phi_I(\Delta)$ is also in $\overline{\calc}_n$, it follows that $J$ does not have a well-defined limit on this set.  It does on the remainder of $\overline{\calAC}_n$, though.

\begin{lemma}\label{shrinking J}  For any $n\geq 3$, $J\co\calAC_n\to\mathbb{R}^+$ extends continuously to $\bigsqcup\phi_I(\calAC_m)$, where $I$ runs over index sets in $\cali_{m,n}$ with $2<m<n$, and for such $I$ that $J\circ\phi_I = J$.  (Here the left-hand ``$J$''  acts on $\calAC_n$ but the right-hand ``$J$'' on $\calAC_m$.)\end{lemma}

\begin{proof}  Since $J\co\calAC_n\to\mathbb{R}$ is symmetric, it suffices to consider sequences in $\calAC_n$ converging to some $(\bd,0,\hdots,0)$ for $\bd\in\calAC_m$, $2<m<n$; i.e. to restrict attention to $\phi_I(\calAC_m)$ where $I = (0,1,2,\hdots,m-1)$.  We will moreover assume that $\bd$ has maximal first entry.  Lemma \ref{reduce} then implies that the equation $d_0 = \ell^n(J,d_1,\hdots,d_{n-1})$ holds on a neighborhood of $\bd$ in $\calAC_m$.  We claim that in fact it holds on a neighborhood of $(\bd,0,\hdots,0)\in\calAC_n$.

The key fact is that one of $d_0 = \ell^n(J,d_1,\hdots,d_{n-1})$ or $d_0 = \ell^s(J,d_1,\hdots,d_{n-1})$ holds at each point of any such neighborhood, where $J = J(d_0,\hdots,d_{n-1})$.  But $\ell^s$ decreases in $J$ on $[D/2,\infty)$, and $\ell^s(D/2)<D$ for $D = \max\{d_i\}_{i=1}^{n-1}$, as calculated in the proof of Lemma \ref{reduce}.  The claim follows.

Lemma \ref{reduce} asserts that $J(\bd) > D/2$, where $D = \max\{d_i\}_{i=1}^{m-1}$.  Since $\ell^n$ increases in $J$ we have $\ell^n(J_0) < d_0 < \ell^n(J_1)$ for any $J_0$, $J_1$ with $D/2 < J_0 < J(\bd) < J_1$, taking $d_0$ here to be the initial entry of $\bd$.  But on inspecting the definition of $\ell^n$, it is clearly a continuous function of $(d_1,\hdots,d_{n-1},J)$ on a neighborhood of $(\bd,0,\hdots,0,J(\bd))$, so the inequality above holds for all $(d_0,\hdots,d_{n-1})\in\calAC_n$ sufficiently near $(\bd,0,\hdots,0)$.  It follows that $J_0 < J(d_0,\hdots,d_{n-1}) < J_1$ at such points.\end{proof}

Below we show that the area function $D_0$ extends continuously to all of $\overline{\calAC}_n$.  Since the radius function $J$ extends continuously to $\calAC_n \cup \bigcup \phi_I(\calAC_m)$ the diagonal functions $\ell_{i,j}$ of Corollary \ref{diagonals} do too, so by its definition in Corollary \ref{area function} $D_0$ does as well.  Proposition \ref{horocyclic defect} further handles the extension to $\calHC_n$, so what we address below is the $\phi_I(\Delta)$.

\begin{lemma}\label{shrinking diagonals}  For any $n\geq 4$ and $0\leq i_0 < i_1 \leq n-1$ let $I = (i_0,i_1)$.  For any $d\geq 0$,
$$ \lim_{\bd\to \phi_I(d,d)} \ell_{i,j}(\bd) = \left\{\begin{array}{ll}
  0 & \mbox{if $i$ and $j$ are cyclically between $i_0$ and $i_1$; or}\\
  d & \mbox{otherwise} \end{array}\right. $$
Here we say $i$ and $j$ are \mbox{\rm cyclically between} $i_0$ \mbox{\rm and} $i_1$ if either $i_0\leq i, j < i_1$ or if each of $i$ and $j$ is less than $i_0-1$ (taken to be $n-1$ if $i_0=0$) or at least $i_1$.\end{lemma}

\begin{proof}  One simply observes that for any cyclic polygon $\{x_0,\hdots,x_{n-1}\}$ with side length collection $(d_0,\hdots,d_{n-1})$ near $\phi_I(d,d)$, there is a sequence of edges joining $x_i$ to $x_j$ such that all but at most one has length near zero; and there is one with length near $d$ if and only if $i$ and $j$ are not cyclically between $i_0$ and $i_1$.  The result thus follows from the triangle inequality.\end{proof}

\begin{corollary}\label{shrinking D}  For any $n\geq 3$, $D_0\co\calAC_n\to\mathbb{R}^+$ extends continuously to $\overline{\calAC}_n$, satisfying $D_0\circ\phi_I = D_0$ for any $I\in\cali_{m,n}$ with $2< m<n$, and $D_0 \equiv 0$ on $\phi_I(\Delta)$ for any $I\in\cali_{2.n}$.\end{corollary}

\begin{proof}  For $A$ as in Lemma \ref{triangle area} it is straightforward to compute:
$$ \lim_{(x,y,z)\to(0,0,0)} A(x,y,z) = 0 = \lim_{(x,y,z)\to(d,d,0)} A(x,y,z) $$
Inspecting the definition of $D_0$ in Proposition \ref{smooth D} one sees that its values on a sequence approaching $\phi_I(d,d)$ are sums of limits of the forms above.
\end{proof}

\bibliographystyle{plain}
\bibliography{cyclic_bib}

\end{document}